\input amstex
\documentstyle{amsppt}
\magnification=\magstep 1
\catcode`\@=11
\def\nologo{\let\logo@=\relax}
\catcode`\@=\active
\nologo
\vsize6.7in

\topmatter
\abstract  
Hochster   established the existence of a commutative noetherian ring $\Cal R$  and a universal resolution $\Bbb U$ of the form $0\to \Cal R^{e}\to \Cal R^{f}\to \Cal R^{g}\to 0$  such that for any commutative noetherian ring $S$ and any resolution  $\Bbb V$ equal to  $0\to S^{e}\to S^{f}\to S^{g}\to 0$,  there exists a unique ring homomorphism $\Cal R\to S$ with $\Bbb V=\Bbb U\otimes_{\Cal R} S$. In the present paper we assume that $f=e+g$ and we find a resolution $\Bbb F$ of $\Cal R$ by free $\Cal P$-modules, where $\Cal P$ is a polynomial ring over the ring of integers. The resolution $\Bbb F$ is not minimal; but it is straightforward, coordinate free, and independent of characteristic. Furthermore, one can use $\Bbb F$ to calculate $\operatorname{Tor}^{\Cal P}_{\bullet}(\Cal R, \Bbb Z)$. If  $e$ and $g$ both at least $5$, then $\operatorname{Tor}^{\Cal P}_{\bullet}(\Cal R, \Bbb Z)$ is not a free abelian group; and therefore, the graded betti numbers in the minimal resolution of $\pmb K\otimes_{\Bbb Z} \Cal R$ by free $\pmb K\otimes_{\Bbb Z} \Cal P$-modules depend on the characteristic of the field $\pmb K$. We record the modules in the minimal $\pmb K\otimes_{\Bbb Z} \Cal P$ resolution of $\pmb K\otimes_{\Bbb Z} \Cal R$ in terms of the modules which appear when one resolves divisors over  the  determinantal ring defined by the $2\times 2$ minors of an $e\times g$ matrix.

 \endabstract
\title The resolution of the universal ring for finite length modules of projective dimension two \endtitle 
  \leftheadtext{Andrew R. Kustin}
\rightheadtext{The resolution of the universal ring}
 \author Andrew R. Kustin\endauthor
 \address
Mathematics Department,
University of South Carolina,
Columbia, SC 29208\endaddress
\email kustin\@math.sc.edu \endemail

\endtopmatter

\document
 \bigpagebreak

\flushpar{\bf   Introduction.}\footnote""{2000 {\it Mathematics Subject Classification.} 13D25.}\footnote""{{\it Key words and phrases.} Buchsbaum-Eisenbud multipliers,  Finite free resolution,  Koszul complex,  Universal resolution.}

\medskip
 
Fix  positive integers $e$, $f$, and $g$, with $r_1\ge 1$ and $r_0\ge 0$, for $r_1$ and $r_0$ defined to be $f-e$ and $g-f+e$, respectively.  Hochster \cite{{12}, Theorem 7.2}  established the existence of a commutative noetherian ring $\Cal R$  and a universal resolution $$\Bbb U\:\quad 0\to \Cal R^{e}\to \Cal R^{f}\to \Cal R^{g}\to 0$$ such that for any commutative noetherian ring $S$ and any resolution  $$\Bbb V\: \quad0\to S^{e}\to S^{f}\to S^{g}\to 0,$$ there exists a unique ring homomorphism $\Cal R\to S$ with $\Bbb V=\Bbb U\otimes_{\Cal R} S$. Hochster showed that the universal ring $\Cal R$ is integrally closed and finitely generated as an  algebra over $\Bbb Z$. Huneke \cite{{13}} identified the generators of $\Cal R$ as a $\Bbb Z$-algebra. These generators correspond to the entries of the two matrices from $\Bbb U$ and the $\binom{g}{r_1}$ multipliers from the factorization theorem  of Buchsbaum and Eisenbud \cite{{6}, Theorem 3.1}.
Bruns \cite{{2}} showed that $\Cal R$ is  factorial. Bruns \cite{{3}} also showed that  universal resolutions exist only for resolutions of length at most two. Heitmann \cite{{11}} used Bruns' approach to universal resolutions in his counterexample to the rigidity conjecture.  
Pragacz and Weyman \cite{{20}} found the relations on the generators of $\Cal R$ and used Hodge algebra techniques to prove that $\pmb K\otimes_{\Bbb Z}\Cal R$ has rational singularities when $\pmb K$ is a field of characteristic zero.  
Tchernev  \cite{{21}} used the theory of Gr\"obner bases to generalize and extend all of the above results with special interest in allowing an arbitrary  base ring $R_0$. 

When $r_1$ is equal to 1, then the universal ring $\Cal R$ is the polynomial ring over $\Bbb Z$ with variables which represent entries of the second matrix from $\Bbb U$ together with variables which represent the Buchsbaum-Eisenbud multipliers. In particular, when $g=r_1=1$, then  the Hilbert-Burch theorem, which classifies all resolutions of the form $$0\to \Cal R^{f-1}\to \Cal R^{f}\to \Cal R^1\to 0,$$ is recovered. When $e=1$ and $r_0=0$, then the universal resolution looks like
$$0\to \Cal R^1\to \Cal R^{f}\to \Cal R^{f-1}\to 0,$$
and
the universal ring $\Cal R$ is defined by the generic Herzog ideal of grade $f$   \cite{{1}}. The minimal resolution of $\Cal R$  is given in \cite{{18}}. 
  
The present paper concerns the universal ring $\Cal R$ when $r_0=0$. In this case,  $f=e+g$ and $\Cal R=\Cal P/\Cal J$, for
$\Cal P$ equal to the polynomial ring $\Bbb Z[\frak b,\{v_{jk}\}, \{x_{ij}\}]$, with    $1\le k\le e$, $1\le j\le f$,
and $1\le i\le g$,
where  $\{\frak b\} \cup\{v_{jk}\}\cup\{x_{ij}\}$  is a list of indeterminates over $\Bbb Z$.
 We give $\Cal J$ in the language of \cite{{21}}. Let
$V$ be the $f\times e$ matrix and $X$ be the $g\times f$ matrix  with
entries  $v_{jk}$ and $x_{ij}$, respectively. For each $$\text{partition  
of $\{1,\dots, f\}$ into $A\cup \bar A$ with $|A|=e$ and
$|\bar A |=g$,}\tag {0.1}$$ let $\nabla_{\bar A,A}$ be the sign of the permutation which arranges the elements of $\bar A,A$ into increasing order, $V(A)$ the submatrix of $V$ consisting of the rows from $A$, and $X(\bar A)$ the submatrix of $X$ consisting of the columns from $\bar A$.   In this notation, the ideal which defines the universal ring $\Cal R$ is   $$\Cal J=I_1(XV)+(\{\det X(\bar A)+\nabla_{\bar A,A} \frak b\det V(A)\mid \text{ $A\cup \bar A$ from ({0.1})}\}).\tag {0.2}$$  
We  produce a resolution $\Bbb F$  of $\Cal R$ by free $\Cal P$-modules. 
Let  $E$, $F$, and $G$ be free $\Cal P$ modules of rank $e$, $f=e+g$, and $g$, respectively; and view the matrices $V$ and $X$ as  homomorphisms of $\Cal P$-modules:
 $$ E@> V >> F@> X>> G.$$  The resolution $\Bbb F$ is given in Definition {2.4}. It is infinite, but straightforward, coordinate free, and independent of characteristic. One can view $\Bbb F$ as the mapping cone of the following map of complexes:
$$\CD \bigoplus\limits_{(a,c,d)} D_aE\otimes D_cG^*\otimes \bigwedge^d(E\otimes G^*)\otimes {\tsize \bigwedge}^{a-c+e}F^*\\@VVV\\\ \bigwedge^{\bullet}(E\otimes G^*)
.\endCD\tag {0.3}$$
The map $X\circ V\:E\to G$ induces a map $E\otimes G^*\to R$ which gives rise to an ordinary Koszul complex $$\dots\to {\tsize \bigwedge}^d(E\otimes G^*)\to {\tsize \bigwedge}^{d-1}(E\otimes G^*)\to \dots\ .\tag {0.4}$$ Also, the map $X\:F\to G$ gives rise to the Koszul complex $$\dots \to S_cG\otimes {\tsize \bigwedge}^bF\to S_{c+1}G\otimes {\tsize \bigwedge}^{b-1}F\to \dots $$ and its dual 
$$\dots \to D_{c+1}G^*\otimes {\tsize \bigwedge}^{b-1}F^*\to  D_cG^*\otimes {\tsize \bigwedge}^bF^*\to \dots;\tag {0.5}$$the map $V^*\:F^*\to E^*$ gives rise to the Koszul complex $$\dots\to S_aE^*\otimes {\tsize \bigwedge}^bF^* \to S_{a+1}E^*\otimes {\tsize \bigwedge}^{b-1}F^*\to \dots $$ and its dual $$\dots\to D_{a+1}E\otimes {\tsize \bigwedge}^{f-b+1}F^*\to D_aE\otimes {\tsize \bigwedge}^{f-b}F^*\to \dots;\tag {0.6}$$and  the identity map on $E^*\otimes G$ gives rise to the Koszul complex 
$$\dots \to S_aE^*\otimes S_cG\otimes {\tsize \bigwedge}^d(E^*\otimes G)\to S_{a+1}E^*\otimes S_{c+1}G\otimes {\tsize \bigwedge}^{d-1}(E^*\otimes G)\to \dots \tag {0.7}$$and its dual 
$$\dots \to D_{a+1}E\otimes D_{c+1}G^*\otimes {\tsize \bigwedge}^{d-1}(E\otimes G^*)\to D_aE\otimes D_cG^*\otimes {\tsize \bigwedge}^d(E\otimes G^*)\to \dots \ .\tag {0.8}$$
The bottom complex of ({0.3}) is the Koszul complex ({0.4}). The differential on the top complex of ({0.3})  involves the maps from ({0.4}), ({0.5}), ({0.6}), and ({0.8}).
The map from the top complex to the bottom complex   involves the minors of $X$ and $\frak b$ times minors of  $V$. The proof that $\Bbb F$ is a resolution of $\Cal R$ uses the acyclicity  lemma and induction on $e$. The case $e=1$ is treated in section 3. In this case, $\Cal R$ is defined by a Herzog ideal and the resolution of $\Cal R$ is already known. In section 5, we exhibit a finite free subcomplex $\Bbb G$ of $\Bbb F$ which has the same homology as $\Bbb F$. We complete the proof that $\Bbb G$ and $\Bbb F$ are acyclic in section 6. On many occasions we filter a given complex over a partially ordered set; the basic procedure is sumarized in section 4.    

All calculations in sections 1 -- 6 are made over the ring of integers $\Bbb Z$; however, it is necassary to work over a field in section 7. 
 One can use $\Bbb F$ to calculate $\operatorname{Tor}^{\Cal P}_{\bullet}(\Cal R, \Bbb Z)$.  If  $e$ and $g$ both at least $5$ and $\pmb K$ is a field, then the Hilbert function of the graded vector space $\operatorname{Tor}^{\Cal P\otimes_{\Bbb Z}\pmb K}_{\bullet}(\Cal R\otimes_{\Bbb Z}\pmb K, \pmb K)$ depends on the characteristic of $\pmb K$; and therefore, the graded betti numbers in the minimal resolution of $\pmb K\otimes_{\Bbb Z} \Cal R$ by free $\pmb K\otimes_{\Bbb Z} \Cal P$-modules also depend on the characteristic of the field $\pmb K$. In particular, $\operatorname{Tor}^{\Cal P}_{\bullet}(\Cal R, \Bbb Z)$ is not a free abelain group.  In Theorem {7.6}, we record the modules in the minimal $\pmb K\otimes_{\Bbb Z} \Cal P$ resolution of $\pmb K\otimes_{\Bbb Z} \Cal R$ in terms of the modules which appear when one resolves divisors over  the  determinantal ring defined by the $2\times 2$ minors of an $e\times g$ matrix. Inspired by Theorem {7.6}, Kustin and Weyman \cite{{19}} used the geometric technique for calculating syzygies in order to give a completely different  calculation of the $\pmb K\otimes_{\Bbb Z} \Cal P$ resolution of $\pmb K\otimes_{\Bbb Z} \Cal R$, when $\pmb K$ is a field of characteristic zero. The resolution in \cite{{19}} is expressed in terms of Schur modules.


\bigpagebreak

\flushpar{\bf 1.\quad  Preliminary results.}

\medskip

 In this paper ``ring'' means commutative noetherian ring with one. If $R$ is a ring and $F$ is an $R$-module, then we let $F^*$ denote $\operatorname{Hom}_R(F,R)$. 

 Let $F$ be 
 a free $R$-module of finite rank. We  make much use of the exterior algebra $\bigwedge^{\bullet}F$   and the divided power algebra $D_{\bullet}F$; we make some use of  the symmetric algebra $S_{\bullet}F$ and the tensor algebra $T_{\bullet}F$. In particular, $\bigwedge^{\bullet}F$ and  $\bigwedge^{\bullet}F^*$ are modules over one another. Indeed, if   $\alpha_i\in \bigwedge^iF^*$ and $b_j\in \bigwedge^jF$, then   $$\tsize \alpha_{i}(b_{j})\in \bigwedge^{j-i}F\quad\text{and}\quad  
b_{j}(\alpha_{i})\in \bigwedge^{i-j}F^*.\nopagebreak\tag {1.1}$$ 
(We view  $\bigwedge^iF$  and $D_iF$ to be meaningful for every integer $i$; in particular, these modules are zero whenever $i$ is negative.)  The exterior and divided power  algebras $A$   come equipped with  co-multiplication
$\Delta\: A\to A\otimes A$.
The following   facts are well known; see \cite{{7}, section 1}, \cite{{8}, Appendix}, and \cite{{14}, section 1}.

\proclaim{Proposition {1.2}} Let $F$ be a free module of rank $f$ over a commutative
noetherian ring $R$ and let $b_r\in \bigwedge^{r}F$, $b'_p\in \bigwedge^{p}F$, and $\alpha_q\in\bigwedge^{q}F^{*}$.
\item{\rm(a)}If $p=f$, then $[b_r(\alpha_q)](b'_p)={\vphantom{E^{E^{E}}_{E_{E}}}} b_r\wedge\alpha_q(b'_p)$.
 \item{\rm(b)}If $X\:F\to G$ is a homomorphism of free $R$-modules and $\gamma_{s+r}\in \bigwedge^{s+r}G^*$, then 
$(\bigwedge^{s}X^{*})\left[\left((\bigwedge^{r}X)(b_{r})\right)(\gamma_{s+r})\right]=
b_{r}\left[\left(\bigwedge^{s+r}X^{*}\right)(\gamma_{s+r})\right]$. \endproclaim

\definition{Notation}Let $m$ be an integer. Each pair of elements $(U,Y)$, with $U\in D_{m}E$ and  $Y\in \bigwedge^{m}G^*$,  gives rise to an element of $\bigwedge^m(E\otimes G^*)$, which we denote by 
$U\bowtie Y$. We now give the definition of $U\bowtie Y$. Consider the composition 
$$ D_mE\otimes  T_mG^*@>\Delta\otimes 1>> T_mE\otimes  T_mG^*  @>\psi >>   {\tsize \bigwedge}^m(E\otimes G^*), $$where 
$ \psi\left((U_1\otimes \dots \otimes U_m)\otimes (Y_1\otimes\dots\otimes Y_m){\vphantom{E^{E^{E}}_{E_{E}}}}\right)= (U_1\otimes Y_1)\wedge \dots \wedge (U_m\otimes Y_m)$,
for $U_i\in E$ and $Y_i\in G^*$.  It is easy to see  that the above composition factors through $D_mE\otimes {\tsize \bigwedge}^mG^*$. Let $U\otimes Y\mapsto U\bowtie Y$ be the resulting map from $D_mE\otimes {\tsize \bigwedge}^mG^*$ to $\bigwedge^m(E\otimes G^*)$. Notice, for example, that if $u\in E$ and $Y_i\in G^*$, then 
$$u^{(m)}\bowtie (Y_1\wedge\dots\wedge Y_m)= (u\otimes Y_1)\wedge \dots \wedge(u\otimes Y_m).$$ The map 
$${\tsize \bigwedge}^mE\otimes D_mG^*\to {\tsize \bigwedge^m}(E\otimes G^*),$$which sends $U\otimes Y$ to $U\bowtie Y$, for $U\in {\tsize \bigwedge}^mE$ and $Y\in  D_mG^*$, is defined in a completely analogous manner.
 \enddefinition

    \definition{Definition}
If $Y\:E\to G$ is a homomorphism of free $R$-modules of finite rank, then let $\check{Y}$ be the  element of $\left(E\otimes G^*\right)^*$ which corresponds to $Y$ under the adjoint isomorphism. In other words, $\check{Y}(\varepsilon\otimes \gamma)=[Y(\varepsilon)](\gamma)$.
In light of ({1.1}), we view $\check{Y}$ as a differential on the exterior algebra $\bigwedge^{\bullet}(E\otimes G^*)$.     
\enddefinition
\remark{Remark}If one thinks of $Y$ as a matrix and takes $\varepsilon$ and $\gamma$ to be basis elements  of $E$ and $G^*$, respectively, then $\check{Y}(\varepsilon\otimes \gamma)$ is the corresponding entry of $Y$.   The differential graded algebra $\left(\bigwedge^{\bullet}(E\otimes G^*),\check{Y}\right)$ is the ``Koszul complex'' associated to the  entries of a matrix representation of $Y$.\endremark

  \proclaim{Lemma {1.3}} Suppose $R$ is a polynomial ring over the ring of integers,   $E$ and $M$ are free $R-$modules, and $\varphi\:D_rE\to M$ is an $R-$module homomorphism. If $\varphi(\varepsilon_1^{(r)})=0$ for all $\varepsilon_1\in E$, then $\varphi$ is identically zero.\endproclaim
\remark{Remarks}  If $R\to \bar R$ is any base change, then $\varphi\otimes 1_{\bar R}$ is also identically zero. On the other hand, the lemma would be false if $R$ were allowed to have torsion. 
Indeed, if $R=\Bbb Z/(2)$, $E=Rx\oplus Ry$ has rank $2$,  and $M$ has rank $1$, then $\varphi\:D_3E\to R$, given by $\varphi(x^{(3)})=\varphi(y^{(3)})=0$ and $\varphi(xy^{(2)})=\varphi(x^{(2)}y)=1$, defines  an $R$-module homomorphism with $\varphi(\varepsilon_1^{(3)})=0$ for all $\varepsilon_1\in E$, but $\varphi$ is not identically zero. 
\endremark
\demo{Proof} Every element of $D_rE$ is a linear combination of elements of  the form $$\varepsilon_1^{(a_1)}\cdots \varepsilon_s^{(a_s)}$$ for some positive integers $s$, and $a_1,\dots, a_s$, with $a_1+\dots+a_s=r$, and elements $\varepsilon_1,\dots,\varepsilon_s$ in $E$.  We show that $D_r(E)\subseteq \ker \varphi$ by induction on $s$. The case $s=1$ is the original hypothesis.
Suppose  that all elements   of the above form  are in $\ker \varphi$ for some  $s$. Fix 
the element $Y=\varepsilon_1^{(a_1)}\cdots \varepsilon_s^{(a_s)}\varepsilon_{s+1}^{(a_{s+1})}$ of $D_rE$.
Let $a=a_s+ a_{s+1}$, and $X$ be the element $\varepsilon_1^{(a_1)}\cdots \varepsilon_{s-1}^{(a_{s-1})}$ of $D_{r-a}E$. The induction hypothesis ensures that 
 for each integer $n$, $X(\varepsilon_s+n\varepsilon_{s+1})^{(a)}$ is in $\ker \varphi$. We see that $X(\varepsilon_s+n\varepsilon_{s+1})^{(a)}$ is equal to the product  $$\bmatrix 1&n&n^2&\dots& n^a\endbmatrix \bmatrix \format\l\\ X\varepsilon_s^{(a)}\varepsilon_{s+1}^{(0)}\\ X\varepsilon_s^{(a-1)}\varepsilon_{s+1}^{(1)}\\ X\varepsilon_s^{(a-2)}\varepsilon_{s+1}^{(2)}\\\phantom{X\varepsilon_s^{(a)}}\vdots\\\vspace{1pt} X\varepsilon_s^{(0)}\varepsilon_{s+1}^{(a)}\endbmatrix.$$
 The row vector  in the above product is a row from a Vandermonde matrix. A   matrix argument produces a non-zero integer $N$ so that $NX\varepsilon_s^{(a-i)}\varepsilon_{s+1}^{(i)}\in \ker \varphi$, for all $i$, with $0\le i\le a$. It follows that $N\varphi(Y)=0$ in the free abelian group M; so, $\varphi(Y)=0$.   \qed
\enddemo

\definition{Convention}If $F$ is a free module of rank $f$, then we orient $F$ by  fixing basis elements
$\omega_F\in\bigwedge^{f} F$ and  $\omega_{F^*}\in\bigwedge^{f} F^*$, which are compatible in the sense that  $\omega_F(\omega_{F^*})=1$. \enddefinition

  \definition{Convention}  For each  statement  ``S'', we define  $$\chi(\text{S})= \cases 1,&\text{if S is true, and}\\ 0,&\text{if S is false.}\endcases$$In particular, $\chi(i=j)$ has the same value as the Kronecker delta $\delta_{ij}$. \enddefinition


\bigpagebreak

\flushpar{\bf 2.\quad  The complex $\Bbb F$.}

\medskip

 \definition{Data {2.1}}Let $R$ be a commutative noetherian ring and let $e$, $f$, and $g$ be positive integers which satisfy $f=e+g$.  The complex $\Bbb F$ is built from data $(\frak b,V,X)$ where $\frak b$ is an element of $R$, and $V$ and $X$ are $R$-module homomorphisms:       $$E@> V >> F@>X>> G,$$ with
   $E$, $F$ and $G$   free $R$-modules  of rank  $e$, $f$, and $g$, respectively.
  For integers $a$,   $c$, and $d$, define  
$$A(a,c,d)=D_aE\otimes D_cG^*\otimes {\tsize \bigwedge}^d(E\otimes G^*)\otimes {\tsize \bigwedge}^{a-c+e}F^* \text{ and } 
B(d)={\tsize \bigwedge}^d(E\otimes G^*)
.$$
\enddefinition

\remark{Remark}If  $V=[v_{jk}]$ and $X=[x_{ij}]$ are matrices, with $1\le k\le e$, $1\le j\le f$, and $1\le i\le g$,  and $R$ is the polynomial ring $R_0[\{\frak b\} \cup\{v_{jk}\}\cup\{x_{ij}\}]$,
where  $\{\frak b\} \cup\{v_{jk}\}\cup\{x_{ij}\}$  is a list of indeterminates over a  commutative noetherian ring $R_0$,
 then we say that the data of {2.1} is generic. 
\endremark

\definition{Grading Convention {2.2}}
Let $\Cal A$ be the additive sub-monoid  of $\Bbb Z^2$ which is generated by $(1,0)$, $(0,1)$, and $(-e,g)$. Notice that  $$\text{$0$ is the only invertible element of $\Cal A$.}\tag {2.3}$$
If the data of {2.1} is generic, then the polynomial ring $R=\bigoplus\limits_{a\in \Cal A}R_a$ is graded by $\Cal A$,
 where the variables $v_{jk}$, $x_{ij}$, and
$\frak b$ of $R$  are elements of $R_{(1,0)}$, $R_{(0,1)}$, and $R_{(-e,g)}$, respectively.   Furthermore, if $I(R)$ is the $R$-submodule $\bigoplus\limits_{a\neq 0}R_a$ of $R$, then condition ({2.3}) ensures that $I(R)$ is an ideal  of $R$ and $R/I(R)\cong R_0$. 
\enddefinition
\remark{Remark}A different approach to the grading of $R$ is given by choosing positive integers $d_x$ and $d_v$ with $gd_x-ed_v$ also positive; for example Tchernev,  takes $(\deg \frak b,d_v,d_x)=(g,g,e+1)$. Then the degrees of 
$v_{jk}$, $x_{ij}$, and
$\frak b$  are $d_v$, $d_x$, and $gd_x-ed_v$, respectively.\endremark

\definition{Definition {2.4}} Let  $(\frak b,V,X)$ be the data of {2.1}.  For each integer $i$, the module $\Bbb F_i$ in the complex    
 $(\Bbb F,\pmb d)$ is 
$$\Bbb F_i= B(i) \oplus \bigoplus\limits_{(a,c,d)} A(a,c,d),$$where the parameters satisfy $i=a+c+d+1$.
The  differential $\pmb d\:\Bbb F_i\to \Bbb F_{i-1}$ is defined 
 as follows.
If $z_d\in B(d)$, then 
 $$\pmb d(z_d )= \check{(X\circ V)}(z_d)\in B(d-1),\ \text{and}$$ if $x=\varepsilon_1^{(a)}\otimes \gamma_1^{(c)}\otimes z_d\otimes \alpha_b\in A(a,c,d)$, with $b=a-c+e$, then $\pmb d(x)$ is equal to 
$$ \left\{\matrix \format\l\\  
 \varepsilon_1^{(a-1)}\otimes \gamma_1^{(c)}\otimes z_d\otimes [V(\varepsilon_1)](\alpha_b)\in A(a-1,c,d)\\\vspace{3pt}
-\varepsilon_1^{(a)}\otimes \gamma_1^{(c-1)}\otimes z_d\otimes X^*(\gamma_1)\wedge\alpha_b \in A(a,c-1,d) \\\vspace{3pt}
+(-1)^{a+c}\varepsilon_1^{(a)}\otimes \gamma_1^{(c)}\otimes \check{(X\circ V)}(z_d)\otimes \alpha_b\in A(a,c,d-1)\\\vspace{3pt}
+(-1)^{a+c}\varepsilon_1^{(a-1)}\otimes \gamma_1^{(c-1)}\otimes (\varepsilon_1\otimes \gamma_1)\wedge z_d\otimes \alpha_b\in A(a-1,c-1,d+1) \\\vspace{3pt}
  +(-1)^{a+d} \chi(c=0) \varepsilon_1^{(a)}\bowtie \left[({\tsize \bigwedge}^{f-b}X)(\alpha_{b}[\omega_F])\right](\omega_{G^*})\wedge z_d\in B(a+d)\\\vspace{3pt}
+(-1)^{d}\chi(a=0) \frak b\cdot   [(\bigwedge^bV^*)(\alpha_b)](\omega_E)  \bowtie \gamma_1^{(c)}\wedge z_d\in B(c+d).
\endmatrix \right.$$
\enddefinition

The map  $$A(a,c,d)\to A(a-1,c,d)$$ is the composition
$$\eightpoint \align D_aE\otimes M\otimes {\tsize \bigwedge}^bF^*@>\Delta\otimes 1\otimes 1>> D_{a-1}E\otimes E\otimes M\otimes {\tsize \bigwedge}^bF^*   @>1\otimes V\otimes 1\otimes 1>> D_{a-1}E\otimes F\otimes M\otimes {\tsize \bigwedge}^bF^*\\
@>>> D_{a-1}E\otimes M\otimes F\otimes {\tsize \bigwedge}^bF^*
 @> 1\otimes ({1.1})\otimes 1 >> D_{a-1}E\otimes M\otimes   {\tsize \bigwedge}^{b-1}F^*,\endalign $$
where $M=D_cG^*\otimes \bigwedge^d(E\otimes G^*)$. The maps from $A(a,c,d)$ to $A(a,c-1,d)$ and $A(a-1,c-1,d+1)$ also involve co-multiplication in a divided power algebra.
 \remark{Remark {2.5}}If the data of {2.1} satisfies the grading convention of {2.2}, then the  complex $\Bbb F$ is homogeneous in the  $\Cal A$-grading, provided   
$$\Bbb F_i= B(i)[-i,-i] \oplus \bigoplus\limits_{(a,c,d)} A(a,c,d)[-a-d,-g-c-d],$$where the sum is taken over all parameters $(a,c,d)$ which satisfy $i=a+c+d+1$.
\endremark

\proclaim{Proposition}The maps and modules $(\Bbb F,\pmb d)$ of Definition {2.4} form a complex\endproclaim
\demo{Proof}In light of Lemma {1.3} it suffices to show that $\pmb d\circ \pmb d(x)=0$ for 
$$x=\varepsilon_1^{(a)}\otimes \gamma_1^{(c)}\otimes z_d\otimes \alpha_b\in A(a,c,d)$$ and $b=a-c+e$. The calculation is routine. We pick out a couple of high points. 
The observation that 
if $Y\:E\to G$, then
$$\check Y(\varepsilon_k\bowtie \gamma_1^{(k)}) =\left( [Y^*(\gamma_1)](\varepsilon_k ){\vphantom{E^{E^{E}}_{E_{E}}}}\right)\bowtie \gamma_1^{(k-1)}$$for all $\varepsilon_k\in \bigwedge^kE$ and $\gamma_1\in G^*$, is the key to seeing that   that the $B(a+d-1)$ component of $\pmb d\circ \pmb d(x)$ is zero, when $c=0$, and that the $B(c+d-1)$ component of $\pmb d\circ \pmb d(x)$ is zero, when $a=0$. In the $B(c+d)$ component of $\pmb d\circ \pmb d(x)$, when $a=1$, we use Proposition {1.2}~(b) and (a) to see that  
 $$\split &[({\tsize \bigwedge}^{b-1}V^*)([V(\varepsilon_1)](\alpha_b))](\omega_E)  \bowtie \gamma_1^{(c)}
=
 \left(\varepsilon_1[({\tsize \bigwedge}^{b}V^*) (\alpha_b)]\right)(\omega_E)  \bowtie \gamma_1^{(c)}\\ 
&{} =\left(\varepsilon_1\wedge [({\tsize \bigwedge}^{b}V^*) (\alpha_b)](\omega_E) \right) \bowtie \gamma_1^{(c)}
  =(\varepsilon_1\otimes \gamma_1)\wedge \left([({\tsize \bigwedge}^{b}V^*) (\alpha_b)](\omega_E)   \bowtie \gamma_1^{(c-1)}\right).\endsplit $$
The same type of argument gives 
$$ \left[({\tsize \bigwedge}^{f-b-1}X)([X^*(\gamma_1)\wedge \alpha_{b}][\omega_F])\right](\omega_{G^*})=
 \gamma_1\wedge \left(({\tsize \bigwedge}^{f-b}X) (\alpha_{b} [\omega_F]){\vphantom{E^{E^{E}}_{E_{E}}}}\right)(\omega_{G^*}),$$which is the key to seeing that the $B(a+d)$ component of $\pmb d\circ \pmb d(x)$ is equal to zero when $c=1$. \qed \enddemo

\definition{Definition {2.6}}Let  $(\frak b,V,X)$ be the data of {2.1}.
Define $\lambda$ to be the element $$(-1)^{eg}[({\tsize \bigwedge}^gX^*)(\omega_{G^*})](\omega_F)+\frak b({\tsize \bigwedge}^eV)(\omega_E)$$ of $\bigwedge^eF$ and define $J$ to be the image of the map $$\bmatrix \lambda& \check{(X\circ V)}\endbmatrix\: {\tsize \bigwedge}^eF^*\oplus (E\otimes G^*)@> >> R.$$ 
\enddefinition

\proclaim{Observation}If $(\Bbb F,\pmb d)$ is the complex of {2.4} and $J$ is the ideal of Definition {2.6}, then the homology $\operatorname{H}_0(\Bbb F)$ is equal to  $R/J$. Furthermore, if  $(\Bbb F,\pmb d)$  is formed using the polynomial ring $\Cal P$ and the data $(\frak b, V, X)$ of $({0.2})$, then 
the homology $\operatorname{H}_0(\Bbb F)$ is equal to the universal ring $\Cal R=\Cal P/\Cal J$.\endproclaim

\demo{Proof} The beginning of $\Bbb F$ is $\Bbb F_1\to \Bbb F_0\to 0$, with $$\Bbb F_1=A(0,0,0)\oplus B(1)={\tsize \bigwedge}^eF^*\oplus (E\otimes G^*)\quad\text{and}\quad \Bbb F_0=B(0)=R.$$ The map $E\otimes G^*\to R$ is $\check{(X\circ V)}$, and the  element $\alpha_e\in \bigwedge^eF^*$ is sent to 
$$\left[({\tsize \bigwedge}^g X)(\alpha_e(\omega_F)\right](\omega_{G^*})+\frak b\left[({\tsize \bigwedge}^e V^*)(\alpha_e)\right](\omega_E)=\lambda(\alpha_e)\in R. $$
The first assertion is established. The homomorphisms $X$ and $V$, of the second assertion,  are represented by matrices. Let  $\alpha_e$ be a basis vector in $\bigwedge^eF^*$. The element $\left[({\tsize \bigwedge}^eV^*)(\alpha_e)\right](\omega_E)$ of $R$ is the determinant  of   the submatrix of $V$ determined by the $e$ rows picked out by $\alpha_e$. The element $\alpha_e(\omega_F)$ of $\bigwedge^gF$ picks out the complementary columns of $X$ with the correct sign, and  $\left[({\tsize \bigwedge}^gX)(\alpha_e(\omega_F)\right](\omega_{G^*})$ is the (signed) determinant of this   submatrix of $X$. 
\qed\enddemo

In Theorem {6.1} we prove that $\Bbb F$ is acyclic whenever the data is generic. However, $\Bbb F$ is far from a minimal resolution. On the other hand, it is possible to isolate the part of $\Bbb F$ in which the splitting    occurs. 
To do this, we partition $\Bbb F$ into strands.  Our definition of the strands  is motivated by Remark {2.5}. 

\definition{Definition {2.7}} Let $(\Bbb F,\pmb d)$ be the complex of Definition {2.4} and let $P$ and $Q$ be integers.
The module $A(a,c,d)$ from $\Bbb F$ is in the strand $S(P,Q)$ if $P=a+d$ and $Q=c+d$. The module $B(d)$   is in $S(P,Q)$ if $P=d$ and $Q=d-g$.
We impose the {\it inverse lexicographic order} on $\Bbb Z\times \Bbb Z$. In other words,     ${(P',Q')\le (P,Q)}$, if $Q'< Q$, or if $Q'=Q$ and  $P'\le P$.  \enddefinition

\proclaim{Observation {2.8}}
\item{\rm(a)} If the strand $S(P,Q)$ is non-zero, then $0\le P$ and $-g\le Q-P\le e$.
\item{\rm(b)} As a module, $\Bbb F=\bigoplus_{(P,Q)}S(P,Q)$.
\item{\rm(c)} The decomposition of {\rm(b)} satisfies   hypothesis {4.1}.
\endproclaim

\demo{Proof} The first assertion holds because if  $A(a,c,d)$ is a non-zero summand of $S(P,Q)$, then $0\le b\le f$
and $b=a-c+e=P-Q+e$. Assertion (b) is obvious. 
If $P=d$ and $Q=d-g$, then $B(d)\subset S(P,Q)$ and $\pmb d(B(d))\subset S(P-1,Q-1)$. If $a+d=P$ and $c+d=Q$, then $A(a,c,d)\subset S(P,Q)$ and 
$$\pmb d(A(a,c,d))\subset \left\{\matrix \format\l\\ S(P-1,Q)\oplus S(P,Q-1)\oplus S(P-1,Q-1)\oplus S(P,Q)\\{}\oplus \chi(c=0) S(P,Q+a-g)\oplus \chi(a=0) S(P+c,Q-g).\endmatrix \right. $$
We have already seen that if $c=0$, then $a-g\le 0$.  
Assertion  (c) has been established. \qed
\enddemo

For each $(P,Q)$ in the poset $\Bbb Z\times \Bbb Z$, let  $(S(P,Q),\pmb\partial)$ be the homogeneous strand of $(\Bbb F,\pmb d)$ which is induced by the direct sum decomposition of (b) as described in {4.1}. Let $(\Bbb F,\pmb\partial)$ be the direct sum of the all of the strands  $(S(P,Q),\pmb\partial)$.

\example{Example {2.9}}Fix integers $P$ and $Q$. If $P-g<Q$, then 
 $S(P,Q)$ is
$$0\to A(P,Q,0)@>\pmb\partial>> A(P-1,Q-1,1)@>\pmb\partial>>\dots @>\pmb\partial>> 
A(P-eg,Q-eg,eg)\to 0.$$
If $Q=P-g$, then $S(P,Q)$ is
$$ 0\to A(P,Q,0)@>\pmb\partial>> A(P-1,Q-1,1)@>\pmb\partial>>\dots @>\pmb\partial>> 
A(P-Q,0,Q)@>\pmb\partial>> B(P)\to 0.$$
In all cases the module $A(a,c,d)$ is in position $a+c+d+1$. \endexample


\bigpagebreak

\flushpar{\bf 3.\quad  The case  $e=1$.}

\medskip

\proclaim{Theorem {3.1}}The complex $(\Bbb F,\pmb d)$ of Definition {2.4} is acyclic when the data is generic and $e=1$.\endproclaim
\demo{Proof}In Definition  {3.3} and Proposition {3.4}, we produce $q\:\Bbb P\to \Bbb P\,'$, which  is  a map of acyclic complexes.  We define  a map of complexes $\varphi\:\Bbb F\to \Bbb P$ in   Proposition {3.5}. It is clear that $q\circ \varphi\:\Bbb F\to \Bbb P\,'$ is surjective. In Proposition {3.7}   we identify the kernel of $q\circ \varphi$ as $M+\pmb d M$.     Lemma {3.9} gives an isomorphism  of complexes $\Theta\:(\Bbb F,\pmb d)\to (\Bbb F,\pmb D)$ which carries $M+\pmb d(M)$ to $M+\pmb D(M)$. This lemma also  shows that $M+\pmb D(M)$ is split exact. It follows that $M+\pmb d(M)=\ker(q\circ \varphi)$ is split exact; and therefore, $\operatorname{H}_i(\Bbb F)=\operatorname{H}_i(\Bbb P\,')$ for all $i$; thus, $\Bbb F$ is acyclic. \qed
\enddemo

 \definition{Data {3.2}}Let $R$ be a commutative noetherian ring,    $g$ be a positive integer, and  $f=g+1$.  The complex $\Bbb P$ is built from data $(\frak b,v,X)$,  where $\frak b$ is an element of $R$,  $X\:F\to G$ is an $R$-module homomorphism,        with
  $F$ and $G$   free $R$-modules  of rank    $f$ and $g$, respectively, and $v$ is an element of $F$. 
\enddefinition
\remark{Remark}If  we think as the data of {3.2} as matrices
$v=[v_{j1}]$ and $X=[x_{ij}]$, with $1\le j\le f$ and $1\le i\le g$, and $R$ is the polynomial ring $R_0[\{\frak b\} \cup\{v_{j1}\}\cup\{x_{ij}\}]$, where
   $\{\frak b\} \cup\{v_{j1}\}\cup\{x_{ij}\}$  is a list of indeterminates over a   ring $R_0$,    then we say that the data of {3.2} is generic. 
\endremark

\definition{Definition  {3.3}} Adopt the data $(\frak b,v,X)$ of {3.2}.   
The complex  $(\Bbb P,\pmb d)$, has modules   $\Bbb P_i= \bigwedge^iG^*\oplus \bigwedge^iF^*\oplus\bigwedge^{i-1}G^* $ and differential  
$$\pmb d_i=\bmatrix X(v)&B_i&(-1)^{i+1}\frak b\\0&-v&A_{i-1}\\0&0&X(v) \endbmatrix ,$$where the maps $ A_i\:\bigwedge^iG^*\to \bigwedge^iF^*$ and $B_i\:\bigwedge^iF^*\to \bigwedge^{i-1}G^*$ are given by   $  A_i(\gamma_i)=\bigwedge^{i}X^*(\gamma_i)$ and  $B_i(\alpha_i)= [(\bigwedge^{f-i}X)(\alpha_i[\omega_F])](\omega_{G^*})$. 
The complex $(\Bbb P\,',\pmb d')$ is the same as the complex $(\Bbb P,\pmb d)$, except  
$$   \Bbb P\,'_1 ={\tsize \bigwedge}^1 G^* \oplus {\tsize \bigwedge}^1 F^*, \quad    \Bbb P\,'_0=R,\quad \pmb d_2'=   \bmatrix  X(v)  &B_2&-\frak b\\0&-v&A_1 \endbmatrix,$$  and $ \pmb d_1'= 
\bmatrix  X(v) &B_1+ \frak b v\endbmatrix$.  
The map of complexes $q\:\Bbb P\to \Bbb P\,'$ be given by $q_i$ is the identity map for $2\le i$, $$q_1=\bmatrix 1&0&0\\0&1&0\endbmatrix,\quad\text{and}\quad q_0=\bmatrix 1&-\frak b\endbmatrix.$$ \enddefinition
\proclaim{Proposition  {3.4}} If the data $(\frak b,v,X)$ of {3.2} is generic, then the complexes $\Bbb P$ and $\Bbb P\,'$ are acyclic. 
\endproclaim
\demo{Proof}One may prove this directly or see \cite{{18}, Prop\. 1.2} or \cite{{17}, Theorem 1.3}. \qed
\enddemo

\proclaim{Proposition {3.5}}Let  $(\Bbb F,\pmb d)$ be the complex of Definition {2.4} constructed from data $(\frak b,V,X)$ with $e=1$, and let $\Bbb P$ be the complex of Definition {3.3} constructed from $(\frak b,v,X)$, where $v=V(1)\in F$. If the function   $\varphi\:\Bbb F\to \Bbb P$  is defined by $$\varphi(\gamma_d)=\bmatrix \gamma_d\\0\\0\endbmatrix\quad\text{and}\quad 
 \varphi(x)= \bmatrix 0\\\chi(0\le a)\chi(0=c)(-1)^a \alpha_b\wedge (\bigwedge^dX^*)(\gamma_d)\\ \chi(0=a)(-1)^c (\bigwedge^b v)(\alpha_b)\cdot \gamma_1^{(c)}\wedge \gamma_d\endbmatrix,$$
for $\gamma_d\in B(d)$ and $x=1^{(a)}\otimes \gamma_1^{(c)}\otimes \gamma_d\otimes \alpha_b\in A(a,c,d)$, then $\varphi$ is a map of complexes.  \endproclaim

\demo{Proof}We pick out one high point of this  calculation.   Fix $a$, $b$, $c$, and $d$ with $i=a+c+d+1$ and $b=a-c+1$. Let $$y=(\pmb d_i\circ \varphi_i-\varphi_{i-1}\circ \pmb d_i)(x)\in \Bbb P_{i-1}={\tsize \bigwedge}^{i-1}G^*\oplus {\tsize \bigwedge}^{i-1}F^*\oplus {\tsize \bigwedge}^{i-2}G^*. $$
When $c=0$, the $\bigwedge^{i-1}G^*$ component of $y$ is zero because
  $$\split   B_i\left( \alpha_b\wedge ({\tsize \bigwedge}^dX^*)(\gamma_d)\right)& 
   = 
\left(({\tsize \bigwedge}^{f-i}X) \left[\left(\alpha_b\wedge ({\tsize \bigwedge}^dX^*)(\gamma_d)\right)(\omega_F)\right]\right)(\omega_{G^*})
\\ &{}= (-1)^{bd} 
\left(({\tsize \bigwedge}^{f-i}X) \left[\left(  ({\tsize \bigwedge}^dX^*)(\gamma_d)\right)(\alpha_b(\omega_F))\right]\right)(\omega_{G^*})\\
&{}= (-1)^{bd} 
\gamma_d\wedge  \left[({\tsize \bigwedge}^{f-b}X)(\alpha_b(\omega_F))\right] (\omega_{G^*})
.
\endsplit  $$The final equality follows from  Proposition {1.2}~(b) and (a).
 \qed \enddemo

\proclaim{Lemma {3.6}}Let  $(\Bbb F,\pmb d)$ be the complex of Definition {2.4} constructed from data $(\frak b,V,X)$ with $e=1$. If  $1\le i$,
then there are submodules $K(0,1,i-1)$ and ${L(0,1,i-1)}$ of $A(0,1,i-1)$ and a homomorphism $s_i\:\bigwedge^iG^*\to A(0,1,i-1)$ such that
\item{\rm(a)} $A(0,1,i-1)=K(0,1,i-1)\oplus L(0,1,i-1)$,
\item{\rm(b)} $K(0,1,i-1)$ is the image of $\pmb\partial\:A(1,2,i-2)\to A(0,1,i-1)$  {\rm(}see Example {2.9}{\rm)}, 
\item{\rm(c)} exterior multiplication carries  $L(0,1,i-1)$ isomorphically onto $\bigwedge^iG^*$, and 
\item{\rm(d)} $s_i$ is a splitting of the exterior multiplication map 
$\mu\: A(0,1,i-1) \to \bigwedge^iG^*$.
 \endproclaim
\demo{Proof}The module $A(0,1,i-1)$ is equal to $G^*\otimes \bigwedge^{i-1}G^*$. The map $\mu$ of (d) is surjective; hence, there exists $s_i\:\bigwedge^iG^*\to A(0,1,i-1)$ with $\mu\circ s_i$ equal to the identity map on $\bigwedge^iG^*$.
Let $L(0,1,i-1)=\operatorname{im} s_i$ and $K(0,1,i-1)=\ker \mu$.
 The decomposition $A(0,1,i-1)=\ker \mu\oplus \operatorname{im} s_i$ gives (a). To complete the proof, recall that the complex $$ \tsize \dots \to  D_2G^*\otimes \bigwedge^{i-2}G^* \to D_1G^*\otimes \bigwedge^{i-1}G^* \to  D_0G^*\otimes\bigwedge^{i}G^*\to 0$$ is split exact for  $1\le i$.    \qed  \enddemo
 \proclaim{Proposition {3.7}} Adopt the notation and hypotheses of Proposition {3.5}. Let $M$ be the   submodule $M= \bigoplus 
A(a,c,d)$ of $\Bbb F$, where the sum is taken over all $3$-tuples $(a,c,d)$, with $1\le a$ and $1\le c$.  
Then $\ker (q\circ \varphi)=M+\pmb d(M)$.
\endproclaim
\demo{Proof} It is clear that $M+\pmb d(M)\subseteq \ker (q\circ \varphi)$.  It is not difficult to see that   $\ker (q\circ \varphi)_i=0$ for $0\le i\le 1$, and $M_i=0$ for $2\le i$. Henceforth, we take $2\le i$.    We next show that $$\Bbb F_i=M_i+\pmb d(M_{i+1})+B(i)+A(i-1,0,0)+L(0,1,i-2), {3.8}$$where $L(0,1,i-2)$ is defined in Lemma {3.6}. First of all, it is easy to see that $$\Bbb F_i= M_i+B(i)+A(0,0,i-1)+A(0,1,i-2)+ \sum\limits_{a=1}^{i-1}A(a,0,i-1-a).$$ Indeed, if $A(a,c,d)$ is a summand of $\Bbb F_i$, but is not a summand of $M_i$, then either $a=0$ or $c=0$. If $a=0$, then either $c=0$ (in which case $(a,c,d)=(0,0,i-1)$) or $c=1$ (in which case $(a,c,d)=(0,1,i-2)$). If $c=0$ and $1\le a$, then $A(a,c,d)$ is a summand of 
$\sum\limits_{a=1}^{i-1}A(a,0,i-1-a)$.
Apply $\pmb d$ to $A(1,1,i-2)$, which is a summand of $M_{i+1}$, 
to see that $$A(0,0,i-1)\subseteq \pmb d(M_{i+1})+M_i+A(0,1,i-2)+A(1,0,i-2).$$
Apply $\pmb d$ to $A(1,2,i-3)$, which is a summand of $M_{i+1}$, to see that 
$$A(0,1,i-2)\subseteq L(0,1,i-2)+\pmb d(M_{i+1})+M_i.$$
If $1\le a\le i-2$, then apply $\pmb d$ to $A(a+1,1,i-2-a)$, which is a summand of $M_{i+1}$,  to see that 
$$A(a,0,i-1-a)\subseteq M_i+\pmb d(M_{i+1})+A(a+1,0,i-2-a).$$

 Let $$P_i=\frac{B(i)\oplus A(i-1,0,0)\oplus L(0,1,i-2)}{[B(i)\oplus A(i-1,0,0)\oplus L(0,1,i-2)]\cap[M_i+\pmb d(M_{i+1})]}. $$ Now that ({3.8}) is established, we  know that $  {\Bbb F_i}/[M_i+\pmb d(M_{i+1})]\cong P_i$. On the other hand, 
 the composition 
$$\split &\Bbb P_i \cong B(i)\oplus A(i-1,0,0)\oplus L(0,1,i-2) @>\operatorname{nat}>> P_i \cong \frac {\Bbb F_i}{M_i+\pmb d(M_{i+1})}@>\varphi_i>> \Bbb P_i\endsplit $$
is an isomorphism, where $\operatorname{nat}$ is the natural quotient map. It follows that $$[B(i)\oplus A(i-1,0,0)\oplus L(0,1,i-2)]\cap[M_i+\pmb d(M_{i+1})]=0,$$ and $\ker(q\circ \varphi)= M+\pmb d(M)$. \qed\enddemo

\proclaim{Lemma {3.9}} Adopt the notation and hypotheses of Proposition {3.7}. There exists a differential $\pmb D$ on $\Bbb F$ and a module automorphism $\Theta$ of  $\Bbb F$, such that 
\item{\rm(a)}$(\Bbb F,\pmb D)$ is a complex,
\item{\rm(b)} $\Theta\:(\Bbb F,\pmb d)\to (\Bbb F,\pmb D)$ is an isomorphism of complexes,
\item{\rm(c)} $\Theta$ acts like the identity map on $M$,
\item{\rm(d)} $\Theta(M+\pmb d(M))=M+\pmb D(M)$, and
\item{\rm(e)} $M+\pmb D(M)$ is split exact.
\endproclaim
\demo{Proof}Recall $s_d$ and $K(0,1,d)$ from Lemma {3.6}. We define $\Theta \:\Bbb F\to \Bbb F$. The map $\Theta$ acts like the identity on   each $B(i)$. If $$x=1^{(a)}\otimes \gamma_1^{(c)}\otimes \gamma_d\otimes \alpha_b\in A(a,c,d),$$  then $\Theta(x)$ is equal to 
$$   \left\{\matrix 
\format\l\\
+1^{(a)}\otimes \gamma_1^{(c)}\otimes \gamma_d\otimes \alpha_b\in A(a,c,d)\\
+\chi(c=0)\chi(1\le d)(-1)^d  1^{(a+d)}\otimes \gamma_1^{(c)}\otimes 1\otimes \alpha_b\wedge \bigwedge^dX^*(\gamma_d)\in A(a+d,c,0) \\
-\chi(a=0)\chi(c=0)  1^{(a)}\otimes s_d( \gamma_d)\otimes v(\alpha_b)\in A(a,c+1,d-1).
\endmatrix \right.  $$Let $\pmb D_i=\Theta_{i-1}\circ \pmb d_i\circ \Theta_i^{-1}$. Assertions (a), (b), and (c) are established. Assertion (d) follows from (b) and (c). If the element $x$, from the above display, is in $M_i$, then a straightforward calculation shows that     $\pmb D_i(x)$ is equal to 
$$ \left\{\matrix \format\l\\  
 +\chi(2\le a) 1^{(a-1)}\otimes \gamma_1^{(c)}\otimes \gamma_d\otimes v(\alpha_b)\in A(a-1,c,d) \\
+\chi(a=1)\chi(c=1) 1^{(a-1)}\otimes (\text{id}-s\circ \mu)(\gamma_1\otimes \gamma_d)\otimes v(\alpha_b)\in K(0,1,d)\\
-[\chi(2\le c)+ \chi(c=1)\chi(1\le d)] 1^{(a)}\otimes \gamma_1^{(c-1)}\otimes \gamma_d\otimes X^*(\gamma_1)\wedge \alpha_b\in A(a,c-1,d) \\
+(-1)^{a+c}1^{(a)}\otimes \gamma_1^{(c)}\otimes [X(v)](\gamma_d)\otimes \alpha_b\in A(a,c,d-1)\\
+(-1)^{a+c}1^{(a-1)}\otimes \gamma_1^{(c-1)}\otimes   \gamma_1\wedge \gamma_d\otimes \alpha_b\in A(a-1,c-1,d+1).
\endmatrix \right.$$
For each $i$, with $2\le i$, let 
$$\Bbb S_i= M_{i} \oplus K(0,1,i-2)\oplus \sum\limits_{(a,d)}A(a,0,d),$$
where the sum is taken over all pairs $(a,d)$ with $0\le a$, $1\le d$, and $a+d+1=i$. Observe that $\pmb D(M_i)\subseteq \Bbb S_{i-1}$. The only term which causes any effort is the  term in $A=A(a-1,c-1,d+1)$. 
If $a-1=0$, then $b=a-c+e$ forces    $c\le 2$. On the other hand, $A$ is zero unless   $1\le c$. If $c=1$, then  $A =A(0,0,d+1)\subseteq \Bbb S_{i-1}$. If $c=2$, then $A=A(0,1,d+1)$, but it is clear that $$(-1)^{a+c}1^{(a-1)}\otimes \gamma_1^{(c-1)}\otimes   \gamma_1\wedge \gamma_d\otimes \alpha_b\in K(0,1,d+1)\subseteq \Bbb S_{i-1}.$$  If $1\le a-1$, then either $A$ is in $M_{i-1}\subseteq \Bbb S_{i-1}$ or $c=1$, in which case, we still have  
$ A= A(a-1,0,d+1)\subseteq \Bbb S_{i-1}$. 
 
Let $\Bbb S=\bigoplus_{2\le i}\Bbb S_i$. 
We have shown that $\pmb D(M_i)\subseteq \Bbb S_{i-1}$; hence ${M+\pmb D(M)\subseteq \Bbb S}$. We complete the proof by showing that  $(\Bbb S,\pmb D)$ is split exact and ${\Bbb S\subseteq M+\pmb D(M)}$. In a manner analogous to Definition {2.7}, we partition $\Bbb S$ into strands $\bigoplus_{(P,Q)}\bar S(P,Q)$, where the sum varies over all pairs $(P,Q)$ with $1\le P$ and $1\le Q\le P+1$. 
For parameters $P$ and $Q$, the  summand $X(a,c,d)$ of $\Bbb S$ is in $\bar S(P,Q)$ if   $P=a+d$ and $Q=c+d$, where 
$ X=A$ for all $(a,c,d)$, except $(0,1,d)$, and $X=K$ for $(a,c,d)=(0,1,d)$.
Observe that every summand of $\Bbb S$ lives in exactly one strand.   
Our calculation of $\pmb D(x)$, for $x\in M_i$, shows that
the decomposition $\Bbb S=\bigoplus \bar S(P,Q)$ satisfies the  hypothesis   of {4.1}.
 Furthermore, the homogeneous strand  $\bar S(P,Q)$ of $\Bbb S$  is 
$$  \eightpoint \align&    0\to A(P,Q,0)\to A(P-1,Q-1,1)\to \cdots \to A(P-Q,0,Q)\to 0,\\\intertext{if $Q\le P$; and}  
 & 0\to A(P,Q,0)\to A(P-1,Q-1,1) \to  \cdots \to K(0,1,Q-1)\to 0,\endalign $$if $Q=P+1$. 
These homogeneous strands $\bar S(P,Q)$ are exact because the complexes 
$$\eightpoint \split &0\to D_{Q}G^*\otimes {\tsize \bigwedge}^0G^* \to D_{Q-1}G^*\otimes {\tsize \bigwedge}^1G^*\to \dots\to D_0G^*\otimes {\tsize \bigwedge}^{Q}G^*\to 0,\ \text{and}\\ & 0\to D_{Q}G^*\otimes {\tsize \bigwedge}^0G^* \to D_{Q-1}G^*\otimes {\tsize \bigwedge}^1G^*\to \dots\to K(0,1,Q-1)\to 0\endsplit $$ are split exact    since $1\le Q$.  We conclude that $\Bbb S\subseteq M+\pmb D(M)$ and that $(\Bbb S,\pmb D)$ is a is split exact complex.
\qed \enddemo
 

\bigpagebreak

\flushpar{\bf 4.\quad  Filtrations.}

\medskip

On numerous occasions we consider a filtration on a complex. We are particularly interested in the associated graded object of the filtration, and for that reason we highlight the ultimate associated graded object, even as we set up the filtration. 
\definition{Notation}Let $(\Bbb E,\pmb d)$ be a complex and $\Pi$ be a partially ordered set. 
Suppose that, as a graded module, $\Bbb E=\bigoplus_{p\in \Pi}\Bbb E^{[p]}$ and that 
$$\matrix  \format\l\\\text{for each fixed $p\in \Pi$, the 
 modules and maps  $(\bigoplus\limits_{p'\le p}\Bbb E^{[p']},\pmb d)$}\\\vspace{-5pt}\text{form a subcomplex of $(\Bbb E,\pmb d)$.} \endmatrix \tag{4.1}$$
For each fixed $p\in \Pi$, let $(\Bbb E^{[p]},\pmb\partial)$ be the quotient complex which is given by the following short exact sequence of complexes:
$$0\to (\bigoplus\limits_{p'< p}\Bbb E^{[p']},\pmb d)\to (\bigoplus\limits_{p'\le p}\Bbb E^{[p']},\pmb d)\to (\Bbb E^{[p]},\pmb\partial)\to 0.$$In particular, the map $\pmb\partial_i \:\Bbb E_i^{[p]}\to \Bbb E_{i-1}^{[p]}$ is equal  to   the composition 
$$\Bbb E_i^{[p]}@>\operatorname{incl}>>\Bbb E_i @>\pmb d_{i}>> \Bbb E_{i-1} @>\operatorname{proj} >> \Bbb E_{i-1}^{[p]},$$ and we refer to 
$\pmb\partial$ as the homogeneous part of $\pmb d$ of degree zero with respect to $\Pi$. We refer to each complex $(\Bbb E^{[p]},\pmb\partial)$ as a homogeneous strand of the original complex $(\Bbb E,\pmb d)$. The graded complex associated to the above filtration of $\Bbb E$ is denoted by  $(\Bbb E,\pmb\partial)$ and is equal to $\bigoplus\limits_{p\in\Pi}(\Bbb E^{[p]},\pmb\partial)$. 
\enddefinition

We apply the filtration technique in three settings.   Proposition {4.2}~(a) is   a quick proof of the well-known fact that if the associated graded complex is exact, then so is the original complex. 
Proposition {4.2}~(b) will be used to split an acyclic summand from a complex.  We can look at one homogeneous strand at a time to determine that $\operatorname{im} \pmb\partial_j$ is a summand of $\Bbb E_{j-1}$. Proposition {4.2}~(b) allows us to conclude that the image of the original map $\pmb d_j$ is also a summand of $\Bbb E_{j-1}$. 

\proclaim{Proposition {4.2}}Let $(\Bbb E, \pmb d)$ be a complex of finitely generated projective $R$-mod\-ules and $\Pi$ be a   partially ordered set.  Suppose that $\Bbb E$ may be decomposed as a direct sum $\bigoplus_{p\in \Pi} \Bbb E^{[p]}$ and that this decomposition satisfies  hypothesis {4.1}.  Fix an integer $j$.  \item{\rm(a)} If $\operatorname{H}_j( \Bbb E,\pmb\partial)=0$,  then $\operatorname{H}_j(\Bbb E,\pmb d)=0$.
\item{\rm(b)} If $\operatorname{H}_j( \Bbb E,\pmb\partial)=0$ and $\operatorname{im} \pmb\partial_j$ is a summand of $\Bbb E_{j-1}$,  then $\operatorname{im} \pmb d_j$ is also a summand of $\Bbb E_{j-1}$.
 \endproclaim

 \demo{Proof}  Let $x$ be a non-zero $j$-cycle   of $\Bbb E$. Consider  $x=\sum x^{[p]}$, with $x^{[p]}\in  \Bbb E^{[p]}$, and  let $$ U(x)=\{p\in \Pi\mid x^{[\pi]}=0\ \text{for all $\pi\in \Pi$ with $p\le \pi$}\}.$$ Let $p_0$ be a maximal element of the support of $x$. It is clear that $\pmb\partial(x^{[p_0]})=0$. It follows that there exists $y\in \Bbb E^{[p_0]}$ with $\pmb\partial(y)=x^{[p_0]}$. We see that $ U(x)\subsetneq  U(x-\pmb d y)$. The proof of (a) is completed by induction.
We prove (b). Let $\Bbb E\,'_{j-1}$ be a direct sum  complement  of $\operatorname{im} \pmb\partial_{j}$ in $\Bbb E_{j-1}$. Assertion (a) may be applied to 
$$\overline{\Bbb E}:\quad \Bbb E_{j+1}@>>> \Bbb E_j @>>> \frac{\Bbb E_{j-1}}{\Bbb E\,'_{j-1}} \to 0.$$
We are given that $(\overline{\Bbb E}, \pmb\partial)$ is exact. We conclude that $(\overline{\Bbb E}, \pmb d)$ is exact. It follows readily that  $\Bbb E_{j-1}=\operatorname{im} \pmb d_{j}\oplus \Bbb E\,'_{j-1}$.  $\qed$\enddemo

In \cite{{16}}, we said  that the complex $\Bbb L$ is {\it splittable} if  $\Bbb L$ is the direct sum of two subcomplexes $\Bbb L'$ and $\Bbb L''$, with $\Bbb L'$   split exact, and the differential on $\Bbb L''$   identically zero.
Suppose that $\Bbb L$ is a complex of projective modules and $\Bbb L$ is bounded in the sense that there exists an integer $N$ with $\Bbb L_i=0$ for all $i<N$.  Under these hypotheses, we proved that   $\Bbb L$ is splittable if and only if $\operatorname{H}_j(\Bbb L)$ is projective for all $j$. 

Our third application of the filtration technique is stated in Observation {4.3}. After the notation is set, then the hypothesis is that various homogeneous strands of the complex $\Bbb E$ have been identified and each of these strands contains a splittable substrand. The conclusion is that, in the original non-homogeneous complex $\Bbb E$, each splittable substrand  may be replaced by its homology, at the expense of complicating the differential.

\proclaim{Observation {4.3}} Let $(\Bbb E,\pmb d)$ be a complex of finitely generated projective $R$-mod\-ules and $\Pi$ be a partially ordered set. Suppose that $\Bbb E$ may be decomposed as a direct sum $\bigoplus_{p\in \Pi} \Bbb E^{[p]}$ and that this decomposition satisfies  hypothesis {4.1}.  Suppose that each module $\Bbb E_i^{[p]}$ decomposes into $\Bbb L_i^{[p]}\oplus \Bbb K_i^{[p]}$. Let $\Bbb L_i=\bigoplus_p\Bbb L_i^{[p]}$, $\Bbb K_i=\bigoplus_p\Bbb K_i^{[p]}$, and $\Bbb L=\bigoplus_i \Bbb L_i $. View $\Bbb L$ as a substrand of $(\Bbb E,\pmb\partial)$. 
  If      $ \Bbb L$ is a splittable complex,  then there exists a split exact subcomplex $(\Bbb N,\pmb d)$ of $(\Bbb E,\pmb d)$ such that $\Bbb N$ is a direct summand of $\Bbb E$ as a module,  and $(\Bbb E /\Bbb N)_i\cong   \operatorname{H}_i(\Bbb L )\oplus \Bbb K_i $. 
\endproclaim
\remark{Remark}  We emphasize that ``view $\Bbb L$ as a substrand of $(\Bbb E,\pmb\partial)$'' means that the differential on $\Bbb L$ is 
$$\Bbb L_i@>\operatorname{incl}>> \Bbb E_i @>\pmb\partial >>\Bbb E_{i-1}@> \operatorname{proj}>> \Bbb L_{i-1}.$$
In practice, for a particular choice of $i$ and $p$, one usually takes either $\Bbb L_i^{[p]}$ or $\Bbb K_i^{[p]}$  to be zero. When this practice is in effect, then $\Bbb L$ is easily seen to be a complex. \endremark

\demo{Proof}
The hypothesis guarantees that   $\Bbb L$ decomposes into the direct sum of two subcomplexes $\Bbb P\oplus \Bbb Q$, where $\Bbb Q$ is split exact and $\Bbb P\cong \operatorname{H}(\Bbb L)$.
   For each $i$, let $\Bbb Q_i$ equal $A_i\oplus B_i$, where $B_i$ is equal to the image of $\Bbb Q_{i+1}$ in $\Bbb L$.   We see that the differential in $\Bbb L$   carries $A_i$ isomorphically onto $B_{i-1}$.  Observe that  $\Bbb E_i=A_i\oplus B_i\oplus \Bbb P_i\oplus   \Bbb K_i$ and that the composition
$$A_i @>\operatorname{incl}>>\Bbb E_i @>\pmb d_{i}>> \Bbb E_{i-1} @>\operatorname{proj} >> B_{i-1}\tag {4.4}$$ is an isomorphism for each $i$. The second assertion holds because  the homogeneous part of ({4.4}) is  an isomorphism. Define $\Bbb N$ to be $\bigoplus_i\Bbb N_i$ and $\Bbb M$ to be $\bigoplus_i\Bbb M_i$, with $\Bbb N_i=A_i+   \pmb d_{i+1}(A_{i+1})$ and $\Bbb M_i=\Bbb P_i\oplus \Bbb K_i$. Use the decomposition $\Bbb E_i=A_i\oplus B_i\oplus \Bbb M_i$ to produce the projection maps
$$\pi_i^{B}\:\Bbb E_i\to B_i\quad\text{and}\quad \pi_i^{\Bbb M}\:\Bbb E_i\to \Bbb M_i.$$ Let $\theta_{i-1}\: B_{i-1}\to A_i$ be the inverse of the map of ({4.4}); $\psi_i\:\Bbb E_i\to \Bbb M_i$ be 
$$\psi_i=\pi_i^{\Bbb M}\circ (1-\pmb d_{i+1}\circ \theta_i\circ \pi_i^{B});$$ and $m_i\:\Bbb M_i\to \Bbb M_{i-1}$ be the composition 
$$\Bbb M_i@>\text{incl}>> \Bbb E_i@>\pmb d_i>>\Bbb E_{i-1}@>\psi_{i-1}>> \Bbb M_{i-1}.$$ A straightforward calculation (see, for example, \cite{{15}, Prop\.~7.2} or \cite{{14}, Prop\.~3.14}) shows that
$$0 \to (\Bbb N,\pmb d|_{\Bbb N})@>\text{incl}>>  (\Bbb E,\pmb d)@> \psi>> (\Bbb M,m)\to 0$$
is short exact sequence  of complexes, and that $\Bbb N$ fulfills all of the requirements.  We notice, for future reference, that the decomposition $\Bbb M=\bigoplus_{p\in \Pi} \Bbb M^{[p]}$ also satisfies   hypothesis {4.1}. \qed \enddemo

\bigpagebreak

\flushpar{\bf 5.\quad  Split a huge summand  from $\Bbb F$.}

\medskip
In Corollary {5.3} we exhibit a finite free subcomplex $\Bbb G$ of $\Bbb F$ which has the same homology as $\Bbb F$.

Fix the complexes $(\Bbb F,\pmb d)$ and $(\Bbb F,\pmb\partial)$ of Definition {2.4} and {2.9}.  We define complexes $(\Bbb P(a_0,c_0,d_0),\pmb\partial)$ and $(\Bbb E,\pmb\partial)$. Each of the new complexes is a quotient of $(\Bbb F,\pmb\partial)$ under the natural quotient map. In particular, $\Bbb P_i$ and $\Bbb E_i$ are  defined for all integers $i$. If we don't specify a value for one of these modules, then the module is automatically equal to zero. The position of the module $A(a,c,d)$ is $a+c+d+1$ in every complex which contains it. 
 Let $$\bar A(a,c,d)=\frac{A(a,c,d)}{\pmb\partial (A(a+1,c+1,d-1))}\quad\text{and}\quad 
\bar B(d)= \frac{B(d)}{\pmb\partial (A(g,0,d-g))}.$$

\definition{Definition} If $a_0$, $c_0$, and $d_0$ are integers, with $a_0$ and $c_0$  non-negative, then let $(\Bbb P(a_0,c_0,d_0),\pmb\partial)$ be the complex
$$0\to A(a_0+d_0,c_0+d_0,0)@>\pmb\partial>> \dots @>\pmb\partial>> A(a_0+1,c_0+1,d_0-1)@>\pmb\partial>> A(a_0,c_0,d_0)\to 0.$$
\enddefinition 
If $P=a_0+d_0$ and $Q=c_0+d_0$, then the complex $\Bbb P(a_0,c_0,d_0)$ is a quotient of the homogeneous strand  $(S(P,Q),\pmb\partial)$ of Observation {2.8}.  The homogeneous strands $(S(P,Q),\pmb\partial)$ have been studied extensively, under a slightly different name, in \cite{{16}}. The exact connection between the two notations is  
$$S(P,Q)=\cases
\Bbb M(P,Q)\otimes \bigwedge^{b}F^*,&\text{if $-e\le P-Q\le g-1$, and}\\
\widetilde{\Bbb M}(P,Q)\otimes \bigwedge^{b}F^*,&\text{if $P=Q+g$},\\
\endcases$$for $b=P-Q+e$. The differential $\pmb\partial$ of $S(P,Q)$ is equal to the tensor product of the differential of $\Bbb M(P,Q)$ or $\widetilde{\Bbb M}(P,Q)$ with the identity map on $\bigwedge^bF^*$. 
The following calculations are Corollaries 5.1 and 5.2 of  \cite{{16}}. 
\proclaim{Theorem {5.1}}
\item{\rm (a)}Assume $1-e\le P-Q\le g-1$. If either  $eg-g+1\le Q$ or $eg-e+1\le P$, then $S(P,Q)$ is split exact.

\item{\rm (b)} If $Q=e+P$, then  $S(P,Q)$  has free homology   which is equal to   $\bar A(0,e,P)$; furthermore, if $eg-e+1\le P$, then $\bar A(0,e,P)=0$.

\item{\rm (c)} If $g+Q=P$, then $S(P,Q)$ has free homology which is equal to $\bar B(P)$; furthermore, if $eg-e+1\le P$, then $\bar B(P)=0$.

\item{\rm (d)} Fix integers $a$, $c$, and $d$. Assume that $1-e\le a-c\le g-1$.  If $g-1\le a$ or $e-1\le c$,  then the complex $\Bbb P(a,c,d)$  has free homology equal to  $\bar A(a,c,d)$.
 \endproclaim
\proclaim{Lemma {5.2}}
If  $(\Bbb E,\pmb\partial) $ is the complex 
 $0\to \bigoplus A(a,c,d)\to 0$,  where the parameters satisfy $eg\le a+c+d$,
 then $\Bbb E$ is a splittable complex and $\operatorname{H}(\Bbb E)$ is equal to the free module  $\bigoplus\bar A(a,c,d)$, where the parameters satisfy $eg=a+c+d$.   \endproclaim

\demo{Proof} Observe that $\Bbb E$ is equal to the direct sum 
$$ \bigoplus\limits_{eg=a_0+c_0+d_0}\Bbb P(a_0,c_0,d_0)\oplus\bigoplus\limits_{ eg+1\le a_0+c_0+d_0\atop{0=a_0c_0}}\Bbb P(a_0,c_0,d_0),$$with the parameters $a_0$, $c_0$, and $d_0$ all non-negative. Let $\Bbb P$ be the strand $ \Bbb P(a_0,c_0,d_0)$ of $\Bbb E$, and let $P=a_0+d_0$.
If $eg+1\le a_0+c_0+d_0$ and $a_0c_0=0$, then Theorem {5.1} yields that $\Bbb P$ is split exact. If $a_0=g+c_0$, then $c_0$ must be zero; hence, $eg+1\le P$ and part (c) applies. If $c_0=a_0+e$, then $a_0$ must be zero; hence, $eg-e+1\le P$ and (b) applies. If $1-g\le c_0-a_0\le e-1$, then $c_0\le e-1$; hence, $eg-e+1\le P$ and (a) applies.  
 If $a_0+c_0+d_0=eg$, then $\Bbb P$ has free homology equal to $\bar A(a_0,c_0,d_0)$. Indeed, 
either (d) applies directly or else (a)--(c) yield that $\Bbb P$ is a truncation of a split exact complex of free modules. In the later case the homology of $\Bbb P(a_0,c_0,d_0)$ is the projective module $\bar A(a_0,c_0,d_0)$. All of the complexes are made over $\Bbb Z$ and transfered to the arbitrary ring $R$ of {2.1} by way of base change. In any event $\bar A(a_0,c_0,d_0)$ is a free module. 
 \qed\enddemo
\proclaim{Corollary {5.3}} Let $(\Bbb F,\pmb d)$ be the complex of Definition {2.4}. If $A(a,c,d)$ is a summand of $\Bbb F$ with $a+c+d=eg$, then there exists a free submodule  $A'(a,c,d)$ of $A(a,c,d)$ so that $$A(a,c,d)=A'(a,c,d)\oplus \pmb\partial A(a+1,c+1,d-1).$$   If $\Bbb G$ is the   submodule 
$$\bigoplus\limits_{a+c+d=eg}A'(a,c,d)\oplus \bigoplus\limits_{a+c+d\le eg-1}A(a,c,d)\oplus \bigoplus_i B(i)$$ of $\Bbb F$, then $\operatorname{H}_i(\Bbb G,\pmb d)=\operatorname{H}_i(\Bbb F,\pmb d)$ for all integers $i$.
\endproclaim

\demo{Proof}Lemma {5.2} guarantees the existence of  $A'(a,c,d)$. 
Apply Proposition {4.2} to see that $\Bbb F_{eg+1}=\operatorname{im} \pmb d_{eg+2}\oplus \bigoplus \limits_{a+c+d=eg}A'(a,c,d)$ and that the cokernel of the inclusion $\Bbb G\subseteq \Bbb F$ is split exact. 
 \qed\enddemo


\bigpagebreak

\flushpar{\bf 6.\quad  The complex $\Bbb F$ is acyclic.}

\medskip

\proclaim{Theorem {6.1}}If the data $(\frak b,V,X)$ of {2.1} is generic, then the complex $(\Bbb F,\pmb d)$ of Definition {2.4}   is acyclic.\endproclaim
\demo{Proof}Corollary {5.3} ensures the existence of a 
subcomplex $\Bbb G$ of $\Bbb F$ such that $\Bbb G$ consists of free modules, $\Bbb G$ has length $eg+1$, and $\Bbb F/\Bbb G$ is split exact.  According to the acyclicity lemma \cite{{7}, Cor\.~4.2}, it suffices to show that $\Bbb G_P$ is acyclic for all prime ideals $P$ of $R$ with $\operatorname{grade} P<eg+1$. Thus, it suffices to show that $\Bbb F_P$ is acyclic for all prime ideals $P$ of $R$ with $\operatorname{grade} P<eg+1$. The ideal $I_1(V)$ has grade $ef\ge eg+1$. If $v$ is an entry of a matrix representation of $V$,  then  Lemma  {6.2}  shows that $\Bbb F_v$ is isomorphic to the complex $\Bbb F$ of Lemma {6.3}, which is  built with generic data over the ring $R_0[v^{-1}]$. The complex $\Bbb F$ of Lemma {6.3} has the same homology as $\Bbb F\,'\otimes \Bbb K$, where   $\Bbb F\,'$ is made from   generic data, with     $e-1$ in place of $e$, and $\Bbb K$ is the Koszul complex associated to the sequence of $g$   new indeterminates. Induction on $e$ completes the result. The base case is Theorem {3.1}. \qed\enddemo 

\proclaim{Lemma {6.2}} 
  Form the complex $(\Bbb F,\pmb d)$ using the data $(\frak b,V,X)$  of {2.1}.  Let $\theta$ and $\tau$ be automorphisms of $F$ and $E$, respectively.  
  Form   $(\Bbb F, \pmb d')$ using   $(u'\frak b,\theta\circ V,   X\circ \theta^{-1})$ and form $(\Bbb F, \pmb d'')$ using $(u''\frak b,V\circ \tau,X)$, where $u'=(({\tsize \bigwedge}^{f}\theta^{-1})[\omega_F])(\omega_{F^*})$ and  $u''$ is $[(\bigwedge^{e}\tau^{-1} )  (\omega_E)](\omega_{E^*})$.   Then the complexes
$(\Bbb F,\pmb d)$, $(\Bbb F,\pmb d')$, and $(\Bbb F,\pmb d'')$  are   isomorphic to one another.\endproclaim

 \demo{Proof}Define  $\varphi'\: (\Bbb F,\pmb d)\to (\Bbb F,\pmb d')$ and  $\varphi''\: (\Bbb F,\pmb d)\to (\Bbb F,\pmb d'')$, by
$\varphi'(  z_d)= u'z_d$ in $B(d)$, $\varphi''( z_d )=({\tsize \bigwedge}^d(\tau^{-1}\otimes 1))(z_d)\in B(d)$, 
$$ \split \varphi'(x)&{}= \varepsilon_1^{(a)}\otimes \gamma_1^{(c)}\otimes z_d\otimes ({\tsize \bigwedge}^b\theta^{*-1})(\alpha_b)\in A(a,c,d), \ \text{and}\\ \varphi''(x)&{}= (D_a\tau^{-1})(\varepsilon_1^{(a)})\otimes \gamma_1^{(c)}\otimes ({\tsize \bigwedge}^d(\tau^{-1}\otimes 1))(z_d)\otimes \alpha_b\in A(a,c,d),\endsplit$$ 
for $z_d\in B(d)$, and 
$x=\varepsilon_1^{(a)}\otimes \gamma_1^{(c)}\otimes z_d\otimes \alpha_b\in A(a,c,d)$. It is not difficult to see that $\varphi'$ and $\varphi''$ both are isomorphisms of complexes. 
  \qed\enddemo

\proclaim{Lemma {6.3}}Let $(\frak b,V,X)$ be the data of {2.1}. Suppose that 
 $E= E'\oplus E''$ and ${F=F'\oplus F''},$ with $E''=R\varepsilon$ and $F''=Rf$. Suppose further that
 $$\gather V=\left[\matrix V'&0\\0&V''\endmatrix\right], \quad\text{and}\quad X=\bmatrix X'&X''\endbmatrix,\quad\text{where} \\ E''@>V''>> F'', \quad E'@>V'>>F'@>X'>> G,\quad\text{and}\quad F''@>X''>>G\endgather $$ are $R$-module homomorphisms, and $V''(\varepsilon)=f$. Form the complexes $(\Bbb F,\pmb d)$ and $(\Bbb F',\pmb d')$ using the data $(\frak b,V,X)$ and   $(\frak b,V',X')$, respectively. Let $x''=X''(f)\in G$. Then there exists a split exact complex $\Bbb L$ and a short exact sequence of complexes:   $$0\to (\Bbb F\,',\pmb d')\otimes ({\tsize \bigwedge}^{\bullet}G^*,x'')  \to \Bbb F\to \Bbb L\to 0.$$
 \endproclaim

\demo{Proof}Take $\omega_E=\varepsilon\wedge\omega_{E'}$ and $\omega_F=f\wedge\omega_{F'}$. Let  $\alpha$ be the element of $F^*$ with $\alpha(F')=0$ and $\alpha(f)=1$,  $A(a_1,c,d_1;a_2,d_2;b_1,b_2)$ equal $$D_{a_1}E'\otimes D_{a_2}E''\otimes D_cG^*\otimes {\tsize \bigwedge}^{d_1}(E'\otimes G^{*})\otimes {\tsize \bigwedge}^{d_2}(E''\otimes G^{*})\otimes {\tsize \bigwedge}^{b_1}F^{\prime*}\otimes {\tsize \bigwedge}^{b_2}F^{\prime\prime*},$$ 
and $ B(d_1;d_2)={\tsize \bigwedge}^{d_1}(E'\otimes G^{*})\otimes {\tsize \bigwedge}^{d_2}(E''\otimes G^{*})$. 
We see that $$A(a,c,d)=\sum  A(a_1,c,d_1;a_2,d_2;b_1,b_2)\quad\text{and}\quad B(d)=\sum B(d_1;d_2),$$where the first sum varies over all  tuples $(a_1,a_2,d_1, d_2,b_1,b_2)$ with $a_1+a_2=a$, $b_1+b_2=a-c+e$, and $d_1+d_2=d$, and the second sum varies over all tuples $(d_1, d_2)$ with $d_1+d_2=d$. 
The complex  $(\Bbb F,d)$ is built using the modules $A(a,c,d)$ and $B(d)$.
The complex $(\Bbb F\,',\pmb d')$ is built using  the modules $$A'(a,c,d)= A(a,c,d;0,0;b_1,0)\quad\text{and}\quad B'(d)=B(d;0),$$ where $b_1=a-c+e-1$.
The differential in  the complex  $\left((\Bbb F\,',\pmb d')\otimes ({\tsize \bigwedge}^{\bullet}G^*,x''),D{\vphantom{E^{E^{E}}_{E_{E}}}} \right)$ is given by $$D(a\otimes b)=\pmb d'(a)\otimes b+(-1)^{|a|+1}a \otimes x''(b).$$ Define  $\varphi: \left({\vphantom{E^{E^{E}}_{E_{E}}}} (\Bbb F\,',\pmb d')\otimes ({\tsize \bigwedge}^{\bullet}G^*,x''),D\right)\to (\Bbb F,\pmb d)$ by 
$$ \varphi \left(z_{d_1}\otimes \gamma_{d_2}\in B'(d_1)\otimes {\tsize \bigwedge}^{d_2}G^*\right)= (-1)^{d_2}z_{d_1}\wedge(\varepsilon^{(d_2)} \bowtie \gamma_{d_2})\in B(d_1;d_2),$$and $\varphi\left(\varepsilon_1^{(a_1)}\otimes \gamma_1^{(c)}\otimes z_{d_1}\otimes \alpha_{b_1}\otimes \gamma_{d_2}\in A'(a_1,c,d_1)\otimes \bigwedge^{d_2}G^* \right)$ is equal to 
$$\varepsilon_1^{(a_1)}\otimes \gamma_1^{(c)}\otimes z_{d_1}\wedge  (\varepsilon^{(d_2)}\bowtie \gamma_{d_2})\otimes \alpha_{b_1}\wedge\alpha\in A(a_1,c,d_1;0,d_2;b_1,1).$$ It is clear that $\varphi$ is injective. We see that  $\varphi$ is a   map of complexes, because  $\pmb d\circ \varphi$ and $\varphi\circ D$ both carry the element 
$z_{d_1}\otimes \gamma_{d_2}$ of $B'(d_1)\otimes {\tsize \bigwedge}^{d_2}G^*$
to 
$$(-1)^{d_2}\check{(X'\circ V')}(z_{d_1})\wedge (\varepsilon^{(d_2)}\bowtie \gamma_{d_2}) 
+(-1)^{d_1+d_2}z_{d_1}\wedge \left(\varepsilon^{(d_2-1)}\bowtie x''(\gamma_{d_2})\right),$$
 and  carry the element 
$\varepsilon_1^{(a_1)}\otimes \gamma_1^{(c)}\otimes z_{d_1}\otimes \alpha_{b_1}\otimes \gamma_{d_2}$ of $A'(a_1,c,d_1)\otimes {\tsize \bigwedge}^{d_2}G^* $ 
to 
$$\left \{\matrix \format\l\\  
\phantom{+} \varepsilon_1^{(a_1-1)}\otimes \gamma_1^{(c)}\otimes z_{d_1}\wedge  (\varepsilon^{(d_2)}\bowtie \gamma_{d_2})\otimes [V'(\varepsilon_1)](\alpha_{b_1})\wedge\alpha
\\
-\varepsilon_1^{(a_1)}\otimes \gamma_1^{(c-1)}\otimes z_{d_1}\wedge  (\varepsilon^{(d_2)}\bowtie \gamma_{d_2})\otimes X^{\prime *}(\gamma_1)\wedge\alpha_{b_1}\wedge\alpha  \\
+(-1)^{a_1+c}\varepsilon_1^{(a_1)}\otimes \gamma_1^{(c)}\otimes \check{(X'\circ V')}(z_{d_1})\wedge  (\varepsilon^{(d_2)}\bowtie \gamma_{d_2})\otimes \alpha_{b_1}\wedge\alpha
 \\
+(-1)^{a_1+c+d_1}\varepsilon_1^{(a_1)}\otimes \gamma_1^{(c)}\otimes z_{d_1}\wedge (\varepsilon^{(d_2-1)}\bowtie x''(\gamma_{d_2}))\otimes \alpha_{b_1}\wedge\alpha
\\
+(-1)^{a_1+c}\varepsilon_1^{(a_1-1)}\otimes \gamma_1^{(c-1)}\otimes (\varepsilon_1\otimes \gamma_1)\wedge z_{d_1}\wedge  (\varepsilon^{(d_2)}\bowtie \gamma_{d_2})\otimes \alpha_{b_1}\wedge\alpha\\
  +(-1)^{a_1+d_1+d_2}\delta_{c0} \varepsilon_1^{(a_1)}\bowtie [(\bigwedge^{f-b_1-1}X')(\alpha_{b_1}[\omega_{F'}])](\omega_{G^*})\wedge z_{d_1}\wedge  (\varepsilon^{(d_2)}\bowtie \gamma_{d_2})\\
+(-1)^{d_1+d_2}\chi(a_1=0) \frak b\cdot   [(\bigwedge^{b_1}V^{\prime *})(\alpha_{b_1})](\omega_{E'})  \bowtie \gamma_1^{(c)}\wedge z_{d_1}\wedge  (\varepsilon^{(d_2)}\bowtie \gamma_{d_2}).
\endmatrix \right.$$
 The cokernel of $\varphi$ is the direct sum  of all $A(a_1,c,d_1;a_2,d_2;b_1,b_2)$ such that either $0<a_2$ and $1=b_2$; or
$0\le a_2$ and $0=b_2$.  
Decompose $\operatorname{coker} \varphi$ into a direct sum of strands $\Cal S(P,Q)$, where $$A(a_1,c,d_1;a_2,d_2;b_1,b_2)\in \Cal S(P,Q)\ \text{if }  P=a_2-b_2+d_2\quad\text{and}\quad  Q=a_1+c+d_1.$$
Impose the inverse lexicographic order of {2.7} on $\{(P,Q)\}$. It is easy to see that 
the decomposition $\operatorname{coker} \varphi=\bigoplus \Cal S(P,Q)$ satisfies hypothesis  {4.1} and that each strand $(\Cal S(P,Q),\pmb\partial)$
is equal to  the split exact sequence 
$$\bigoplus 0\to A(a_1,c,d_1;a_2,d_2;b_1,1)@>\cong>> A(a_1,c,d_1;a_2-1,d_2;b_1,0)\to 0,$$where the sum varies over all tuples $(a_1,c,d_1;a_2,d_2;b_1,1)$ with $a_2+d_2=P+1$, $a_1+c+d_1=Q$, and $1\le a_2$. Apply Proposition {4.2}. \qed
\enddemo
 
The ideal $\Cal J$ of ({0.2}) is generically perfect of grade $eg+1$.  The  theorem about the transfer of perfection (see, for example, \cite{{5}, Theorem 3.5})
 tells us that $\Bbb F\otimes_{\Cal P}R$ and $\Bbb G\otimes_{\Cal P}R$ are resolutions for any ring $R$ for which $\Cal JR$ is a proper ideal of grade at least $eg+1$.


\bigpagebreak

\flushpar{\bf 7.\quad  The minimal resolution.}

\medskip

\definition {Data {7.1}} Let $\pmb K$ be a field and $\Cal R=\Cal P/\Cal J$ be the universal ring of ({0.2}). 
We write $\Cal P'$ for the polynomial ring $\pmb K\otimes_{\Bbb Z}\Cal P= \pmb K [\frak b,\{v_{jk}\}, \{x_{ij}\}]$, $\Cal J'$ for the image of $\Cal J$  in $\Cal P'$, and $\Cal R'$ for $\pmb K\otimes_{\Bbb Z} \Cal R=\Cal P'/\Cal J'$. Let $E_0$, $F_0$, and $G_0$ be vector spaces of dimension $e$, $f$, and $g$ over $\pmb K$, and $E=E_0\otimes_{\pmb K}\Cal P'$, $F=F_0\otimes_{\pmb K}\Cal P'$, and $G=G_0\otimes_{\pmb K}\Cal P'$ be the corresponding free $\Cal P'$-modules. 
Let $(\Bbb F,\pmb d)$ be the complex of {2.4}  built using the data $(\frak b, V,X):$
$$E@>V >> F@>X >> G,$$
where  $V=[v_{jk}]$ and $X=[x_{ij}]$ are matrices.  \enddefinition 

In Theorem {7.6}, we record the modules of the minimal $\Cal A$-homogeneous resolution of $\Cal R'$ by free $\Cal P'$-modules.   
There are two steps in our proof of Theorem {7.6}. In the first step, Lemma {7.2},  we apply the technique of Observation {4.3} to the present situation. 
The other step is the calculation of the homology of the homogeneous strands of $\Bbb F$. This step was largely carried out in \cite{{16}}. Most of the modules that comprise the resolution of {7.6} are equal to  modules which arise when one resolves divisors  of a determinantal ring defined by the $2\times 2$ minors of an $e\times g$ matrix.

\proclaim{Lemma {7.2}}Adopt the hypotheses of {7.1}. If  $\Bbb X$ is the minimal $\Cal A$-homogeneous resolution of $\Cal R'$ as a $\Cal P'$-module, then $\Bbb X$  and $\operatorname{H}(\Bbb F\otimes_{\Cal P'}\pmb K)\otimes_{\pmb K}\Cal P'$
are isomorphic  as  $\Cal A$-graded $\Cal P'$-modules. 
\endproclaim

\demo{Proof}The homology of  $\Bbb F\otimes_{\Cal P'}\pmb K$ is free over $\pmb K$, since $\pmb K$ is a field; and therefore,  $\Bbb F\otimes_{\Cal P'}\pmb K$ is equal to the direct sum of two graded subcomplexes, one of which is split exact and the other has zero differential. Use the graded version of Nakayama's Lemma to pull this decomposition back to $\Bbb F$. At this point, the summand $\Bbb F_i$ of $\Bbb F$ has been decomposed as the direct sum $A_i\oplus B_i\oplus C_i$  of graded free $\Cal P'$-modules and the graded differential $\pmb d_i\:\Bbb F_i \to \Bbb F_{i-1}$, which looks like $$\bmatrix \pmb d_i^{11}&\pmb d_i^{12}&\pmb d_i^{13}\\\pmb d_i^{21}&\pmb d_i^{22}&\pmb d_i^{23}\\\pmb d_i^{31}&\pmb d_i^{32}&\pmb d_i^{33}
\endbmatrix,$$ has $\pmb d^{13}_i\otimes_{\Cal P'} \pmb K\:C_i\otimes_{\Cal P'} \pmb K\to A_{i-1}\otimes_{\Cal P'} \pmb K$ is an isomorphism and $\pmb d^{k\ell}_i\otimes_{\Cal P'} \pmb K=0$ if $(k,\ell)\neq (1,3)$. It follows that the map $\pmb d^{13}_i\:C_i\to A_{i-1}$, of graded free $\Cal P'$-modules, is also an isomorphism. Let $\Bbb C$ be the subcomplex  $\bigoplus_iC_i+ \pmb d(\bigoplus_iC_i)$ of $\Bbb F$. 
Ordinary row and column operations produce a 
short exact sequence of complexes of graded free $\Cal P'$-modules
$$0\to (\Bbb C,\pmb d)@>\operatorname{incl}>> (\Bbb F,\pmb d)\to (\Bbb X,x)\to 0,$$
where $\Bbb X=\bigoplus_i B_i$,   $x_i\:B_i\to B_{i-1}$ is $\pmb d_i^{22}-\pmb d_i^{23}\circ (\pmb d_{i}^{13})^{-1}\circ \pmb d_i^{12}$, and 
 $(\Bbb F,\pmb d)\to (\Bbb X,x)$ induces an isomorphism on homology.
 The explicit form of the map $x_i$ guarantees that the differential in $\Bbb X\otimes_{\Cal P'} \pmb K$ is zero.  \qed\enddemo

Now we must identify the homology of $\Bbb F\otimes_{\Cal P'}\pmb K$. We begin by recalling  the bi-graded structure on $\operatorname{Tor}$. 
\definition{Definition} If $\frak P=\bigoplus_i \frak P_i$ is a graded ring,   and $A=\bigoplus_i A_i$ and $B=\bigoplus_i B_i$ are graded $\frak P$-modules, then the module $\operatorname{Tor}^{\frak P}_{\bullet}(A,B)$ is a bi-graded $\frak P$-module. Indeed, if $$\Bbb Y\: \dots \to Y_1 \to Y_0\to A$$ is a $\frak P$-free  resolution of $A$,  homogeneous of degree zero, then 
$$\operatorname{Tor}_{p,q}^\frak P(A,B)=\frac{\ker [(Y_p\otimes B)_q\to (Y_{p-1}\otimes B)_q]}{\operatorname{im} [(Y_{p+1}\otimes B)_q\to (Y_{p}\otimes B)_q]}.$$\enddefinition

\definition{Notation {7.3}}Adopt the notation of {7.1}. Define the $\pmb K$-vector spaces $$\matrix \format\r&\ \l&\ \l\\ \Cal N(a,c,d)&{}=S_aE_0^*\otimes S_cG_0\otimes {\tsize \bigwedge}^d(E_0^*\otimes G_0),\\
\Cal M(a,c,d)&{}=D_aE_0\otimes D_cG_0^*\otimes {\tsize \bigwedge}^d(E_0\otimes G_0^*),&\text{and}\\
B_0(i)&{}=\bigwedge^i(E_0\otimes_{\pmb K}G_0^*).\endmatrix$$ 
   The identity map on $E_0^*\otimes G_0$ induces Koszul complexes of the form  $$\cdots \to \Cal N(a-1,c-1,d+1)\to \Cal N(a,c,d)\to 
\Cal N(a+1,c+1,d-1)\to \cdots \ . \tag {0.7}$$ The  $R$-dual of ({0.7}) is  $$\cdots \to 
\Cal M(a+1,c+1,d-1)\to\Cal M(a,c,d)\to \Cal M(a-1,c-1,d+1)\to\cdots\ .\tag {0.8}$$
Fix integers $P$ and $Q$. Let $\Bbb N(P,Q)$ and $\Bbb M(P,Q)$ be the above complexes  when $a+d=P$ and $c+d=Q$; that is,  $\Bbb N(P,Q)$ is 
$$0\to \Cal N(P-eg,Q-eg,eg)@>>> \dots @>>> \Cal N(P-1,Q-1,1)@>>> \Cal N(P,Q,0)\to 0,$$
and $\Bbb M(P,Q)$ is 
$$0\to \Cal M(P,Q,0)@>>> \Cal M(P-1,Q-1,1)@>>> \dots@>>>  \Cal M(P-eg,Q-eg,eg)\to 0.$$
 If $P=g+Q$, then let $\widetilde{\Bbb M}(g+Q,Q)$ be the augmented complex 
$$0\to \Cal M(g+Q,Q,0)@>>>\dots @>>>\Cal M(g,0,Q)@>\gamma >> B_0(g+Q),$$where $\gamma(U\otimes 1\otimes Z)=(U\bowtie \omega_{G^*})\wedge Z$.   
\enddefinition
\remark{Remark} In \cite{{16}}, the homology of ({0.7}) at $\Cal N(a,c,d)$  is denoted by $\operatorname{H}_{\Cal N}(a,c,d)$ and   the cohomology of ({0.8}) at $\Cal M(a,c,d)$ is called  
 $\operatorname{H}_{\Cal M}(a,c,d)$. \endremark

Recall the strands $(S(P,Q),\pmb\partial)$ of $\Bbb F$, which were introduced at the end of section 2. Observe that 
$$S(P,Q)\otimes_{\Cal P'}\pmb K=\cases \Bbb M(P,Q)\otimes_{\pmb K} {\tsize \bigwedge}^{P-Q+e}F_0^*&\text{if $P<Q+g$}\\
\widetilde{\Bbb M}(P,Q)\otimes_{\pmb K}{\tsize \bigwedge}^{f}F_0^* 
&\text{if $P=Q+g$}.\endcases\tag {7.4}$$

Most of the modules that comprise the resolution of {7.6} are equal to  modules which arise when one resolves divisors  of a determinantal ring defined by $2\times 2$ minors.   Let $S$ be the ring $S^{\pmb K}_{\bullet}E_0^*\otimes_{\pmb K} S^{\pmb K}_{\bullet}G_0$, $T$ be the subring $$T=\sum_m S_mE_0^*\otimes S_mG_0$$ of $S$, and for each integer $\ell$, let $M_{\ell}$ be the $T$-submodule  
$$M_{\ell}=\sum\limits_{m-n=\ell} S_mE_0^*\otimes S_nG_0$$ of $S$. Give $S$ a grading by saying that $S_mE_0^*\otimes S_nG_0$ has grade $n$, for all $m$ and $n$. We see that  $T$ is a graded ring, and  $\bigoplus M_{\ell}$ is a direct sum decomposition of $S$ into graded $T$-submodules. In particular, $$\text{the graded summand of degree $n$ in $M_{\ell}$ is $S_{n+\ell}E_0^*\otimes S_nG_0$.}$$  Let $\frak P$ be the polynomial ring $S_{\bullet}^{\pmb K}(E_0^*\otimes _{\pmb K}G_0)$. The ring  $\frak P$ is   graded; each element of $S_n(E_0^*\otimes G_0)$ is homogeneous of grade $n$.  The identity  map on $E_0^*\otimes G_0$ induces a graded ring homomorphism  from $\frak P$ onto $T$. Each graded  $T$-module is automatically  a graded $\frak P$-module. Notice that $\frak P\otimes_{\pmb K}\bigwedge_{\pmb K}^{\bullet}(E_0^*\otimes G_0)$ is a homogeneous resolution of $\pmb K$ by free $\frak P$-modules; and therefore,  
 $\operatorname{Tor}^{\frak P}_{\bullet}(M_{\ell},\pmb K)$ is the homology of  ${M_{\ell}\otimes_{\pmb K}\bigwedge_{\pmb K}^{\bullet}(E_0^*\otimes G_0)}$; indeed,  $$\operatorname{Tor}^{\frak P}_{p,q}(M_{\ell},\pmb K)=\operatorname{H}_{\Cal N}(\ell+q-p,q-p,p).$$The homology $\operatorname{Tor}^{\frak P}_{p,q}(M_{\ell},\pmb K)$ is a   $\pmb K$-vector space  for all $p$, $q$, and $\ell$. It follows that  $$\operatorname{H}_{\Cal M}(a,c,d)=\operatorname{H}_{\Cal N}(a,c,d)=\operatorname{Tor}_{d,c+d}^{\frak P}(M_{a-c},\pmb K)\tag {7.5}$$ for all integers $a$, $c$, and $d$.
One may view $\frak P$ as a polynomial ring over $\pmb K$ in $eg$ indeterminates.  The ring $T$ is the determinantal ring defined by the $2\times 2$ minors of the $e\times g$ matrix of indeterminates. The   divisor class group of $T$ is $\Bbb Z$ and  $\ell\mapsto [M_{\ell}]$ is an isomorphism from $\Bbb Z$ to $\operatorname{C}\!\ell\,(T)$. Much more information about the modules $M_{\ell}$ may be found in \cite{{4}}.

We are ready to record the modules in the minimal $\Cal P'$-resolution of $\Cal R'$. Almost all of these modules appear in the minimal $\frak P$-resolution of $M_{\ell}$ for some $\ell$, with $-e\le \ell\le g$. 
The other type of module that appears in {7.6} is the cokernel of the map $\gamma$ of ({7.3}), let $\bar B_0(i)=\frac {B_0(i)}{\operatorname{im} \Cal M(g,0,i-g)}$.  Write $\pmb \alpha$ for $(e-1)(g-1)$.
 
\proclaim{Theorem {7.6}} Adopt the hypotheses of {7.1}.    Then the minimal $\Cal A$-homogeneous resolution 
$$\Bbb X\:\quad 0\to \Bbb X_{eg+1} \to \dots \to \Bbb X_0\to \Cal R'\to 0$$ of $\Cal R'$ by free $\Cal P'$-modules has   $\Bbb X_i$ equal to 
$$\tsize \left\{\matrix\format\l\\  \chi(i\le eg-e)\Cal P'\otimes_{\pmb K} \bar B_0(i)[-i,-i]\\\vspace{5pt}\phantom{\chi(i\le eg-e)}\oplus\\\vspace{5pt}\bigoplus\limits_{(p,q,\ell)} \Cal P'\otimes_{\pmb K}\left[\operatorname{Tor}_{p,q}^{\frak P}(M_{\ell},{\pmb K})\otimes_{\pmb K}\bigwedge^{\ell+e}F_0^*\right][-\ell-q,-g-q]\\\vspace{5pt}\phantom{\chi(i\le eg-e)}\oplus\\\vspace{5pt} \chi(i\le eg+1) \Cal P'\otimes_{\pmb K}\operatorname{Tor}_{i-1-e,i-1}^{\frak P}(M_{-e},{\pmb K})[e+1-i,1-i-g]
,\endmatrix \right.$$where 
the sum is taken over all parameters $(p,q,\ell)$ with 
$$1-e\le \ell\le g-1, \quad  0\le p\le \pmb \alpha,$$$$ p+\operatorname{max}\{-\ell,0\}\le q\le p+\min\{g-1-\ell,e-1\}\quad\text{and}\quad 
i=\ell+2q-p+1.$$
\endproclaim

\demo{Proof}We apply Lemma {7.2} to the complex $(\Bbb F,\pmb d)$ of Theorem {6.1}. Observation {2.8} shows that  $\Bbb F\otimes_{\Cal P'}\pmb K$ splits into the following direct sum of complexes:
$$\bigoplus_{(P,Q)} (S(P,Q)\otimes_{\Cal P'} \pmb K,\pmb\partial),$$where the sum is taken over all integers $(P,Q)$ with $0\le P$ and $-g\le Q-P\le e$. Apply ({7.4}) to see that 
$$\Bbb F\otimes_{\Cal P'}\pmb K=\left\{\matrix \format\l \\
\bigoplus\limits_{P}\left[ \widetilde{\Bbb M}(P,P-g)\right][-P,-P]\\\vspace{5pt} \phantom{\bigoplus\limits_{P}\left[ \widetilde{\Bbb M}(P,P-g)\right]}\oplus \\\vspace{5pt}\bigoplus\limits_{(P,Q)}\left[ \Bbb M(P,Q)\otimes_{\pmb K}{\tsize \bigwedge}^{P-Q+e}F^*_0\right][-P,-g-Q],\endmatrix\right. \tag {7.7}$$ 
where the top sum is taken over  all integers $P$, and the bottom sum is taken over all pairs $(P,Q)$ with $-e\le P-Q\le g-1$. It is shown in \cite{{16}} that the homology of each complex $\widetilde{\Bbb M}(P,P-g)$ is free and is concentrated in the position of $B(P)$. It follows immediately that the contribution of the top line of ({7.7}) to $\operatorname{H}_i(\Bbb F\otimes_{\Cal P'} \pmb K)$  is  $\bar B_0(i)[-i,-i]$. 
The  contribution of the bottom line of ({7.7}) to  $\operatorname{H}_i(\Bbb F\otimes_{\Cal P'} \pmb K)$  is 
$$\bigoplus\limits_{(P,Q)}\left[ \operatorname{H}_{\Cal M}(a,c,d)\otimes_{\pmb K}{\tsize \bigwedge}^{P-Q+e}F^*_0\right][-P,-g-Q],$$where $(P,Q)$ continue to satisfy 
$-e\le P-Q\le g-1$, and 
the parameters $(a,c,d)$ satisfy $a+d=P$, $c+d=Q$, and $a+c+d+1=i$. 
Apply ({7.5}) and reparameterize by letting $\ell=P-Q$, $q=Q$, and $p=d$ to see that  
$$\tsize \Bbb X_i=\left\{\matrix\format\l\\  \Cal P'\otimes_{\pmb K} \bar B_0(i)[-i,-i]\\\vspace{5pt}\phantom{\chi(i\le eg-e)}\oplus\\\vspace{5pt}\bigoplus\limits_{(p,q,\ell)} \Cal P'\otimes_{\pmb K}\left[\operatorname{Tor}_{p,q}^{\frak P}(M_{\ell},{\pmb K})\otimes_{\pmb K}\bigwedge^{\ell+e}F_0^*\right][-\ell-q,-g-q],\endmatrix \right.$$where 
the sum is taken over all parameters $(p,q,\ell)$ with 
$$i=\ell+2q-p+1\quad\text{and}\quad -e\le \ell\le g-1.$$

 We know from \cite{{16}}, that 
$$0\to \operatorname{H}_{\Cal M}(g,0,i-g)\to B_0(i)\to \operatorname{H}_{\Cal N}(0,e,eg-e-i)\to 0\tag {7.8}$$
is an exact sequence. Moreover, $$\bar B_0(i)=\operatorname{coker}\left(\operatorname{H}_{\Cal M}(g,0,i-g)\to B_0(i)\right)=\operatorname{H}_{\Cal N}(0,e,eg-e-i).$$The module $\operatorname{H}_{\Cal N}(0,e,eg-e-i)$ is zero if $eg-e-i<0$; and therefore,
$$\bar B_0(i)=\chi(i\le eg-e)\bar B_0(i).$$

Fix $\ell=-e$. The module $\operatorname{Tor}_{p,q}^{\frak P}(M_{-e},{\pmb K})$ is equal to $\operatorname{H}_{\Cal M}(q-p-e,q-p,p)$ and \cite{{16}} tells us that this module is zero unless $q=p+e$. Furthermore, if $q=p+e$, then there is a short exact sequence,
$$0\to \operatorname{H}_{\Cal M}(0,e,p)\to B_0(e+p)\to \operatorname{H}_{\Cal N}(g,0,eg-e-g-p)\to 0.$$The module $B_0(e+p)$ is zero if $eg<e+p$.  If $i=\ell+2q-p+1$ and $\operatorname{Tor}_{p,q}^{\frak P}(M_{-e},{\pmb K})\neq 0$, then 
$$q=i-1,\quad i-e-1=p, \quad\text{and}\quad i\le eg+1.$$When $\ell=-e$, the contribution of 
$$\bigoplus\limits_{(p,q,\ell)} \Cal P'\otimes_{\pmb K}\left[\operatorname{Tor}_{p,q}^{\frak P}(M_{\ell},{\pmb K})\otimes_{\pmb K} {\tsize \bigwedge}^{\ell+e}F_0^*\right][-\ell-q,-g-q]$$ to $\Bbb X_i$ is
$$\chi(i\le eg+1)\Cal P'\otimes_{\pmb K}\operatorname{Tor}_{i-e-1,i-1}^{\frak P}(M_{-e},{\pmb K})[e+1-i,-g-i+1].$$

 Fix $\ell$ with  $1-e\le \ell\le g-1$. The module $$\operatorname{Tor}_{p,q}^{\frak P}(M_{\ell},\pmb K)=\operatorname{H}_{\Cal M}(q-p+\ell,q-p,p)$$ and \cite{{16}, Thm\. 1.1} tells us that this module is isomorphic to 
$$\operatorname{H}_{\Cal N}(g-1-(q-p+\ell),e-1-(q-p),\pmb \alpha-p).$$
It follows that $\operatorname{Tor}_{p,q}^{\frak P}(M_{\ell},{\pmb K})$ is zero unless
$$0\le q-p+\ell\le g-1,\quad  0\le q-p\le e-1,\quad\text{and}\quad 0\le p\le \pmb \alpha; $$ and therefore, $\operatorname{Tor}_{p,q}^{\frak P}(M_{\ell},{\pmb K})$ is zero unless
$$p+\operatorname{max}\{-\ell,0\}\le q\le p+\min\{g-1-\ell,e-1\}\quad\text{and}\quad 0\le p\le \pmb \alpha.
\qed$$
\enddemo

We record some explicit versions of the resolution $\Bbb X$ of Theorem {7.6}. Keep in mind that if $p+p'=\pmb \alpha-1$, then \cite{{16}, Theorem 2.1} shows that 
the dimension of 
$$\bar B_0(g+p') \cong \operatorname{Tor}_{p,p+e}^{\frak P}(M_{-e},\pmb K)$$ can be computed from either of the split exact sequences:
$$\matrix \format\l\\ \vspace{5pt} 0\to \operatorname{Tor}_{p,p+e}^{\frak P}(M_{-e},\pmb K)\to \Cal N(0,e,p)\to \dots \to \Cal N(p,p+e,0)\to 0\quad\text{or} \\
0\to \Cal M(g+p',p',0)\to \dots\to \Cal M(g,0,p')\to B_0(g+p') \to \bar B_0(g+p') \to 0.\endmatrix \tag {7.9} $$ If $e=2$, then the Eagon-Northcott and Buchsbaum-Rim complexes (see \cite{{9}, Theorem A2.10} or  \cite{{16}, Section 4}) give the resolution of $M_{\ell}$ for $-1\le \ell$.
In particular, if $-1\le \ell$, then
$$\dim_{\pmb K}\operatorname{Tor}_{p,q}^{\frak P} (M_{\ell},\pmb K)=\cases \dim S_{\ell-p}E_0\otimes \bigwedge^pG_0&\text{if $q=p$ and $p\le \ell$}\\\dim D_{p-\ell-1}E_0^*\otimes \bigwedge^{p+1}G_0&\text{if $q=p+1$ and $\ell+1\le p$}\\0&\text{otherwise}.\endcases$$
\example{Example {7.10}} To economize space, we write  $$T(p,q,\ell) \text{ for }\Cal P'\otimes_{\pmb K}\left[\operatorname{Tor}_{p,q}^{\frak P}(M_{\ell},{\pmb K})\otimes_{\pmb K}{\tsize \bigwedge}^{\ell+e}F_0^*\right].$$ If $e=g=2$, then $\Bbb X$ is 
$$\alignat{8}  
&\text{module}&&\quad\Bbb X_i&&\quad\text{twist}&&\quad\text{rank}&&\quad\quad\text{module}&&\quad\Bbb X_i &&\quad\text{twist}&&\quad\text{rank}\\ 
&\Cal P'\otimes_{\pmb K}\bar B_0(0)&&\quad 0&&\quad [0,0]&&\quad  1 &&\quad\quad T(2,4,-2)&&\quad 5&&\quad [-2,-6]&&\quad 1 \\ 
&\Cal P'\otimes_{\pmb K}\bar B_0(1)&&\quad 1&&\quad [-1,-1]&&\quad  4&&\quad\quad T(1,3,-2)&&\quad 4&&\quad [-1,-5]&&\quad 4  \\ 
&T(0,0,0)&&\quad 1&&\quad [0,-2]&&\quad 6 &&\quad\quad T(1,2,0)&&\quad 4&&\quad [-2,-4]&&\quad 6 \\
&\Cal P'\otimes_{\pmb K}\bar B_0(2)&&\quad 2&&\quad [-2,-2]&&\quad  3&&\quad\quad  T(0,2,-2)&&\quad 3 &&\quad [0,-4]&&\quad 3 \\ 
&T(0,1,-1)&&\quad 2&&\quad [0,-3]&&\quad 8&&\quad\quad T(1,1,1)&&\quad 3&&\quad [-2,-3]&&\quad  8\\
&T(0,0,1)&&\quad 2&&\quad [-1,-2]&&\quad 8&&\quad\quad     T(1,2,-1)&&\quad 3&&\quad [-1,-4]&&\quad 8.\\
\endalignat$$ 
In other words,  $\Bbb X$ is
$$0\to \smallmatrix \Cal P'[-2,-6]\endsmallmatrix \to \smallmatrix \Cal P'[-2,-4]^6\\\oplus\\ \Cal P'[-1,-5]^4\endsmallmatrix\to \smallmatrix \Cal P'[-2,-3]^8\\\oplus\\\Cal P'[-1,-4]^8\\\oplus\\\Cal P'[0,-4]^3\endsmallmatrix \to \smallmatrix \Cal P'[-2,-2]^3\\\oplus\\\Cal P' [-1,-2]^8\\\oplus\\ \Cal P'[0,-3]^8\endsmallmatrix \to \smallmatrix \Cal P'[-1,-1]^4\\\oplus\\\Cal P'[0,-2]^6\endsmallmatrix \to \smallmatrix \Cal P'. \endsmallmatrix $$
\endexample

\example{Example} Assume that 
 $3\le e$ and $3\le g$.
 The beginning  of $\Bbb X$ is   
$$\eightpoint \alignat{4}  
&\text{module}&&\quad\Bbb X_i&&\quad\text{twist}&&\quad\text{rank}\\ 
&\Cal P'\otimes_{\pmb K}\bar B_0(0)&&\quad 0&&\quad [0,0]&&\quad  1 \\ 
&\Cal P'\otimes_{\pmb K}\bar B_0(1)&&\quad 1&&\quad [-1,-1]&&\quad  eg \allowdisplaybreak\\ 
&T(0,0,0)&&\quad 1&&\quad [0,-g]&&\quad  \binom{f}e \allowdisplaybreak\\ 
&\Cal P'\otimes_{\pmb K}\bar B_0(2)&&\quad 2&&\quad [-2,-2]&&\quad  \binom{eg}2 \allowdisplaybreak\\ 
&T(0,0,1)&&\quad 2&&\quad [-1,-g]&&\quad  e\binom f{e+1} \allowdisplaybreak\\ 
&T(0,1,-1)&&\quad 2&&\quad [0,-(g+1)]&&\quad  g\binom f{e-1} \allowdisplaybreak\\ 
&    \Cal P'\otimes_{\pmb K}\bar B_0(3)&&\quad 3&&\quad [-3,-3]&&\quad  \text{see {7.9}}  \allowdisplaybreak\\
&T(0,0,2)&&\quad 3&&\quad [-2,-g]&&\quad \binom {e+1}2\binom{f}{e+2}     \allowdisplaybreak\\
&T(1,1,1)&&\quad 3&&\quad [-2,-(g+1)]&&\quad  |\Bbb M(2,1)|\binom f{e+1}   \allowdisplaybreak\\
&T(0,2,-2)&&\quad 3&&\quad [0,-(g+2)]&&\quad  \binom f{e-2}\binom {g+1}2    \allowdisplaybreak\\
&T(1,2,-1)&&\quad 3&&\quad [-1,-(g+2)]&&\quad |\Bbb M(1,2)|\binom f{e-1}.    \allowdisplaybreak\\
\endalignat$$ The homology of the complex $\Bbb M(1,2)$ (see {7.3}) is concentrated in one position and is equal to $\operatorname{Tor}_{1,2}^{\frak P}(M_{-1},\pmb K)$; so, the dimension of $\operatorname{Tor}_{1,2}^{\frak P}(M_{-1},\pmb K)$ is equal to the absolute value of the Euler characteristic of $\Bbb M(1,2)$. 

Use ({7.8}) to calculate the contribution of the bottom summand to $\Bbb X_{eg+1}$:
$$\operatorname{Tor}^{\frak P}_{eg-e,eg}(M_{-e},\pmb K)=\operatorname{H}_{\Cal N}(0,e,eg-e)=\pmb K. $$
 The largest index $i$ for which the middle summand contributes to $\Bbb X_i$ is $i=eg$, and the contribution is $\bigwedge^gF^*[-e(g-1),-eg]$ because when $p=\pmb \alpha$, $q=(e-1)g$, and $\ell=g-e$, then 
$$\operatorname{Tor}^{\frak P}_{p,q}(M_{\ell},\pmb K)=\operatorname{H}_{\Cal N}(g-1,e-1,\pmb \alpha)=\operatorname{H}_{\Cal M}(0,0,0)=\pmb K.$$
 The end of $\Bbb X$ is
$$\eightpoint \alignat{4}  
&\text{module}&&\quad\Bbb X_i&&\quad\text{twist}&&\quad\text{rank}\\
&T(eg-e,eg,-e)&&\quad eg+1&&\quad [-(eg-e),-(eg+g)]&&\quad  1 \\ 
\vspace{5pt}\allowdisplaybreak
&T(eg-e-1,eg-1,-e)&&\quad eg&&\quad [-(eg-e-1),-(eg+g-1)]&&\quad  eg \\ 
\vspace{5pt}\allowdisplaybreak
&T(\pmb \alpha,(e-1)g,g-e)&&\quad eg&&\quad [-(eg-e),-eg]&&\quad \tsize \binom{f}e \\
\allowdisplaybreak 
\vspace{5pt} 
&T(eg-e-2,eg-2,-e)&&\quad eg-1&&\quad [-(eg-e-2),-(eg+g-2)]&&\quad  \binom{eg}2 \\ 
\allowdisplaybreak
\vspace{5pt}
&T(\pmb \alpha,(e-1)g,g-e-1)&&\quad eg-1&&\quad [-(eg-e-1),-eg]&&\quad\tsize  e\binom f{e+1} \\
\allowdisplaybreak
\vspace{5pt} 
&T(\pmb \alpha,eg-g-1,g-e+1)&&\quad eg-1&&\quad [-(eg-e),-(eg-1)]&&\quad  \tsize g\binom f{e-1} \\
\allowdisplaybreak
\vspace{5pt}
& T(eg-e-3,eg-3,-e)&&\quad eg-2&&\quad [-(eg-e-3),-(eg+g-3)]&&\quad \operatorname{dim} \bar B_0(3)  \\
\allowdisplaybreak\vspace{5pt}
&T(\pmb \alpha,(e-1)g,g-e-2)&&\quad eg-2&&\quad [-(eg-e-2),-eg]&&\quad \binom {e+1}2\binom{f}{e+2}  \\
\allowdisplaybreak
\vspace{5pt}
&T(\pmb \alpha-1,eg-g-1,g-e-1)&&\quad eg-2&&\quad [-(eg-e-2),-(eg-1)]&&\quad  |\Bbb M(2,1)|\binom f{e+1} \\
\allowdisplaybreak
\vspace{5pt}
&T(\pmb \alpha,eg-g-2,g-e+2)&&\quad eg-2&&\quad   [-(eg-e),-(eg-2)]&&\quad  \tsize\binom f{e-2}\binom {g+1}2   \\
\allowdisplaybreak\vspace{5pt}
&T(\pmb \alpha-1,eg-g-2,g-e+1)&&\quad eg-2&&\quad [-(eg-e-1),-(eg-2)]&&\quad |\Bbb M(1,2)|\binom f{e-1}.    \\
\endalignat$$
\endexample

\example{Example} Let $e=3$ and $g=4$. The module $T(3,4,-1)$ contributes the summand ${\Cal P'[-3,-8]^{210}}$ to $\Bbb X_5$ and the module $T(2,4,-1)$ contributes the summand ${\Cal P'[-3,-8]^{420}}$ to $\Bbb X_6$. By duality, $T(3,4,2)$ contributes the summand ${\Cal P'[-6,-8]^{210}}$ to $\Bbb X_8$ and the module $T(4,4,2)$ contributes the summand ${\Cal P'[-6,-8]^{420}}$ to $\Bbb X_7$.
  These summands can not be predicted if one only knows the Hilbert function of $\Cal R'$. In this particular example, every other summand of the $\Bbb X$ can be correctly predicted from knowledge of the Hilbert function of $\Cal R'$, together with the assumption that the minimal resolution of $\Cal R'$ is as simple as possible.
\endexample

\example{Example} If $\pmb K$ has characteristic  zero, then the $\frak P$-resolution of each module $M_{\ell}$ is known; and therefore, all of the modules in $\Bbb X$ are known in terms of Schur modules; see \cite{{19}}.  The paper \cite{{19}} was inspired by the  present paper; however, the proof is completely different. It uses the geometric method of finding syzygies and is valid only in characteristic zero. The resolution of Example {7.10} may also be found in \cite{{19}}. 
 \endexample

\example{Example} If $e$ and $g$ are both at least $5$, then Hashimoto \cite{{10}} proved that the dimension of  $\operatorname{Tor}^{\frak P}_{3,5}(M_0,\pmb K)$ depends on the characteristic of $\pmb K$; and therefore, the graded betti number $\beta_{8}(5,g+5)$ in $\Bbb X$ depends on the characteristic of $\pmb K$.

\endexample


\medskip

\flushpar {\bf Acknowledgment.}
Thank you to Alexandre Tchernev for getting me started on this project.

\enddocument